\documentclass[times,final]{article}
\usepackage[numbers]{natbib}

\usepackage{framed,multirow}

\usepackage[english]{babel}
\usepackage{lipsum}

\usepackage{amssymb}
\usepackage{latexsym}

\usepackage{url}
\usepackage{xcolor}
\definecolor{newcolor}{rgb}{.8,.349,.1}

\usepackage{lmodern}
\usepackage[T1]{fontenc}
\usepackage{setspace}
\usepackage{xcolor}

\usepackage{pgfplots}
\pgfplotsset{compat=1.11}
\usepackage{tikz,tkz-euclide}
\usetikzlibrary{arrows,calc,patterns,spy,backgrounds}

\usepackage{rotating}
\usepackage[pdfborder={0 0 0}]{hyperref}

\usepackage{amscd}

\usepackage{amsmath,amsthm,mathtools,bm,cancel,relsize}
\usepackage{bbm}
\usepackage{sobolev}
\usepackage{graphics,psfrag,subfigure}
\usepackage{empheq}
\usepackage{colortbl}
\usepackage{booktabs}
\usepackage{enumitem}
\usepackage[margin=1.0in]{geometry}

\newtheorem{thm}{Theorem}[section]

\newtheorem{remark}[thm]{Remark}

\numberwithin{equation}{section}

\usepackage{amsopn}

\newcommand{\Rd}{\color{red}}
\newcommand{\Bd}{\color{blue}}

\newcommand{\Xhath}[1]{\,\hspace{-1pt}\widehat{X}^{#1}_h\hspace{-1pt}\,}
\newcommand{\Xhathb}[1]{\,\hspace{-1pt}\widehat{X}^{#1}_{h,0}\hspace{-1pt}\,}

\newcommand{\XhathT}[1]{\,\hspace{-1pt}\widehat{\tilde{X^{#1}_h}}\hspace{-1pt}\,}
\newcommand{\XhT}[1]{\,\hspace{-1pt}{\widetilde{X}}^{#1}_h\hspace{-1pt}\,}

\let\hat\widehat
\let\tilde\widetilde

\newcommand{\convergenceslope}[6] 
{
\draw [color=black,mark=none,#6]
(axis cs: #1, #3) -- (axis cs: #2, 
{exp((#4*(-1)) * ln((#2)/(#1)) + ln(#3))}) 
-- node [midway, #5, below]{1} (axis cs: #1, 
{exp((#4*(-1)) * ln((#2)/(#1)) + ln(#3))}) 
-- node [midway, #5, left]{${#4}$} (axis cs: #1, #3);
}
\newcommand{\convergenceslopeinv}[6] 
{
\draw [color=black,mark=none,#6]
(axis cs: #1, #3) -- (axis cs: #2, 
{exp((#4*(-1)) * ln((#2)/(#1)) + ln(#3))}) 
-- node [midway, #5, right]{${#4}$} (axis cs: #2, 
#3)
-- node [midway, #5, above]{1} (axis cs: #1, #3);
}
\newcommand{\convergenceslopeh}[6] 
{
\draw [color=black,mark=none,#6]
(axis cs: #2, #3) -- (axis cs: #1, 
{exp((#4*(-1)) * ln((#2)/(#1)) + ln(#3))}) 
-- node [midway, #5, below]{1} (axis cs: #2, 
{exp((#4*(-1)) * ln((#2)/(#1)) + ln(#3))}) 
-- node [midway, #5, right]{${#4}$} (axis cs: #2, #3);
}
\newcommand{\convergenceslopeinvh}[6] 
{
\draw [color=black,mark=none,#6]
(axis cs: #2, #3) -- (axis cs: #1, 
{exp((#4*(-1)) * ln((#2)/(#1)) + ln(#3))}) 
-- node [midway, #5, left]{${#4}$} (axis cs: #1, 
#3)
-- node [midway, #5, above]{1} (axis cs: #2, #3);
}

\def\SC#1{\ifnum `#1<`1 {\Gd \hat {X}^0_h} \else \ifnum  `#1<`2
  \Bd{\hat X^1_h} \else \ifnum  `#1<`3 \Rd{\hat X^2_h} \else \Md{\hat X^3_h} \fi  \fi  \fi}
\def\XC#1{\ifnum `#1<`1 {\Gd {X}^0_h} \else \ifnum  `#1<`2
  \Bd{ X^1_h} \else \ifnum  `#1<`3 \Rd{ X^2_h} \else \Md{ X^3_h} \fi  \fi  \fi}
\def\XCt#1{\ifnum `#1<`1 {\Gd \tilde {X}^0_h} \else \ifnum  `#1<`2
  \Bd{ \tilde X^1_h} \else \ifnum  `#1<`3 \Rd{ \tilde X^2_h} \else \Md{ \tilde X^3_h} \fi  \fi  \fi}

\graphicspath{{./}{./figures/}}

\title{
    High order geometric methods with splines: \\fast solution with explicit time-stepping for Maxwell equations}

\author{Bernard Kapidani$^*$ \and Rafael V\'azquez$^{* \dag}$}
\date{$*$ {\footnotesize Institute of Mathematics, Ecole Polytechnique F\'ed\'erale de Lausanne, Station 8, 1015 Lausanne, Switzerland}\\
$\dag$ {\footnotesize Istituto di Matematica Applicata e Tecnologie Informatiche ``E.\ Magenes'' del CNR, via Ferrata 5, 27100 Pavia, Italy}\\
{\footnotesize (\href{mailto:bernard.kapidani@epfl.ch}{bernard.kapidani@epfl.ch}, \href{mailto:rafael.vazquez@epfl.ch}{rafael.vazquez@epfl.ch})}\\
\vspace{5mm}\today}
\begin{document}

\maketitle



\begin{abstract}
We introduce a high-order spline geometric approach for the initial boundary value problem for Maxwell's equations. The method is geometric in the sense that it discretizes in structure preserving fashion the two de Rham sequences of differential forms involved in the formulation of the continuous system. Both the Amp\`ere--Maxwell and the Faraday equations are required to hold strongly, while to make the system solvable two discrete Hodge star operators are used. By exploiting the properties of the chosen spline spaces and concepts from exterior calculus, a non-standard explicit in time formulation is introduced, based on the solution of linear systems with matrices presenting Kronecker product structure, rather than mass matrices as in the standard literature. These matrices arise from the application of the exterior (wedge) product in the discrete setting, and they present Kronecker product structure independently of the geometry of the domain or the material parameters. The resulting scheme preserves the desirable energy conservation properties of the known approaches. The computational advantages of the newly proposed scheme are studied both through a complexity analysis and through numerical experiments in three dimensions.
\end{abstract}

\section{Introduction} \label{sec:intro}
We are interested in the numerical solution of the initial boundary value problem for the Maxwell equations on a bounded space-time domain $\Omega\times [0, T]$, with $\Omega\subset\mathbb{R}^3$. {The} expression of Maxwell equations in the language of exterior calculus, as compared to vector calculus, allows for a clear separation of topological and geometrical properties of the equations, while also giving a more neat and concise expression of pivotal properties of numerical schemes{. For this reason} we write the system of Maxwell equations in terms of differential forms, which reads 
\begin{align} 
& \partial_t \mathrm{D} =   d \mathrm{H} - \mathrm{J}, \label{eq:ampere}\\
& \partial_t \mathrm{B} = - d \mathrm{E}, \label{eq:faraday}\\
& d\mathrm{D} = 0, \label{eq:divgaugeD} \\
& d\mathrm{B} = 0, \label{eq:divgaugeB}
\end{align}
where we assume for simplicity the absence of free electric charges. The forms $\mathrm{E}$ and $\mathrm{H}$ are respectively the electric and the magnetic field, given as differential 1--forms. $\mathrm{D}$ and $\mathrm{B}$ are respectively the electric displacement and the magnetic induction, given as differential 2--forms, $\mathrm{J}$ is a 2--form representing the electric current, which is a source term of the equations, and $d$ is the exterior derivative, for which a more precise definition will be given in Section~\ref{sec:theory}. Equations \eqref{eq:ampere}--\eqref{eq:divgaugeB} have to be completed with constitutive laws, that take the form
\begin{equation} \label{eq:cont_hodge1}
\begin{aligned}
& \mathrm{B} = \star^1_\mu \mathrm{H},
& \mathrm{D} = \star^1_\epsilon \mathrm{E},
\end{aligned}
\end{equation}
where $\mu$ is the magnetic permeability and $\epsilon$ the electric permittivity, and the $\star^1$ symbol represents the Hodge star operator which maps 1--forms into 2--forms. Again, a more precise definition will be given in Section~\ref{sec:theory}. 

When integrating the electrodynamics \eqref{eq:ampere}--\eqref{eq:faraday} numerically in time, the standard algorithm used in real-life applications is still the finite differences in time domain (FDTD) or its integral variations such as the finite integration technique (FIT)~\cite{ClThWe99,Clemens-Weiland}, or the cell method~\cite{tonti_finite_2001,m_marrone_computational_2001,codecasa_explicit_2008,codecasa_novel_2018}. An alternative framework for Delaunay triangulations is given also by discrete exterior calculus~\cite{Hirani-darcy,hirani}.
These \emph{geometric} methods are first order in space and time and preserve at the discrete level the conservation properties of the continuous Maxwell system. Despite the proliferation of arbitrary order finite element methods (FEM) suitable for the solution of Maxwell equations, achieved by properly discretizing the de Rham complex of differential forms~\cite{HIP02a,AFW06,AFW-2}, their implementation is usually not competitive in terms of computational efficiency. This is due to the amenability to massive parallelization of FDTD structured grids and the simplicity inherent in their implementation, while high order FEM requires the solution of a linear system associated to the mass matrix at each time step. To replicate the success of geometric methods on unstructured meshes, discontinuous Galerkin methods (DG) have become popular. They use arbitrary order polynomials to discretize the electric and magnetic field in each finite element and only weakly enforce their tangential continuity across element boundaries, through numerical fluxes, leading to the solution of linear systems with block diagonal matrices, with one block per element. There are nuances in the way the fluxes are constructed \cite{hesthavenNodalDiscontinuousGalerkin2008}, but in general they are either dissipative fluxes, which sacrifice electromagnetic energy conservation properties, or conservative fluxes, which introduce spurious numerical solutions which cannot be in general eliminated or even recognized as such by the user. High order extensions of geometric methods have been analysed for structured grids in \cite{chungConvergenceSuperconvergenceStaggered2013} and for unstructured ones by one of the authors ~\cite{kapidani_time-domain_2020, kapidani_arbitrary-order_2021}, by recasting them as DG methods on staggered grids.


The language of exterior calculus can also be extended to the emerging framework of isogeometric analysis (IGA). A discrete de Rham complex was first constructed and analysed for tensor-product B-splines in \cite{Buffa_Sangalli_Vazquez,BRSV11} 
and it has been applied in the Galerkin framework for the discretization of Maxwell's equations \cite{ratnani,Corno20161} also in the context of plasma physics \cite{kraus_kormann_morrison_sonnendrucker_2017,Holderied_2021}, and for the development of pointwise divergence free methods for incompressible fluid flow \cite{Buffa_deFalco_Sangalli,EvHu12,EvHu12-2,EvHu12-3,Van17}. The first attempts to fully exploit the concepts of exterior calculus and differential forms for B-splines were made in \cite{back2012} with a dual staggered grid, and in \cite{Hiemstra20141444} where the de Rham complex of B-splines was adapted to the framework of mimetic discretizations as described by Bochev and Hyman \cite{bochevPrinciplesMimeticDiscretizations2006}. Very recently, the isogeometric de Rham complex was combined with a DG approach between conforming patches for the approximation in complex multipatch domains \cite{gucluBrokenFEECFramework2022}.

Building on the listed works, a recent paper by the present authors introduced a new method based on isogeometric differential forms and the construction of two dual de Rham complexes \cite{kapidaniHighOrderGeometric2022}. Compared to other geometric methods, a dual mesh is not explicitly built, and the dual complex is simply defined by a change in polynomial degree, with the same construction introduced in \cite{hiemstra_isogeometric_2011} and applied in \cite{Buffa_2020aa,kapidani2021treecotree} for stable mortar coupling between non-conforming meshes. Thanks to the high continuity of splines, the exterior derivative is rigorously defined in both sequences and given by incidence matrices of a Cartesian grid \cite{ratnani,BSV14}, and the dimension of pairing spaces from the two sequences is always equal. The difference with respect to the approach of \cite{Hiemstra20141444} and \cite{gucluBrokenFEECFramework2022} is that, instead of introducing a discrete version of the co-derivative, which is the adjoint operator of the exterior derivative, we discretize the Hodge star operators \cite{HIP99c} which relate the spaces of the primal and the dual complexes. It was proven in \cite{kapidaniHighOrderGeometric2022} that, when applied to elliptic problems, the method attains high order of convergence, and we also presented numerical evidence showing that the method is spurious free when applied to the Maxwell eigenvalue problem.

The preceding article was a starting point in exploring the general framework for high order geometric methods based on splines, and thus studies the general properties which underpin the approach. The present paper shifts its focus to tailoring the approach to the solution of hyperbolic systems of equations.
The main novelty in the following is to show how to combine the structure from \cite{kapidaniHighOrderGeometric2022} with methods for Kronecker product matrices, similar to the ones applied to preconditioners in \cite{sangalliIsogeometricPreconditionersBased2016}, to achieve a high-order, geometric, and explicit in time approximation of the Maxwell system. In particular, the final matrices for which a linear system has to be solved at each time step will have a Kronecker product structure, independently of the {material parameters or the geometry} of the spatial domain $\Omega$. To solve the associated system we will rely on computing the LU factorization of Gram matrices between univariate spline spaces, which can be then quickly computed, cheaply stored and applied.
Importantly, this differentiates the present work from similar endeavors based on solving linear systems with metric dependent mass matrices, as it is instead done for instance in \cite{ratnani,barhamMimeticDiscretizationMacroscopic2022a}, or the non-conforming approach in \cite{gucluBrokenFEECFramework2022} with mass matrix based blocks.

The outline of the paper is the following. In Section~\ref{sec:theory} we introduce the mathematical tools necessary to fully understand the proposed numerical scheme, mainly revolving around spline complexes of differential forms. In Section~\ref{sec:hodge} we discuss how different choices of discrete Hodge star operators lead to discretization schemes with equivalent conservation properties but notably different algebraic structure. We present the schemes focusing on their matrix form, and comment on their equivalent weak formulations. In Section~\ref{sec:tensorized} we focus on how to exploit the structure of the discrete spaces to efficiently solve the linear systems involved in the explicit time stepping. In Section~\ref{sec:num} we numerically validate high order approximation properties and the efficiency of the scheme. Some summarizing remarks conclude the paper in Section~\ref{sec:fin}.

\section{Preliminaries}\label{sec:theory}

In the present section we introduce the necessary notation to write Maxwell equations in terms of differential forms, including the definition of the Hodge star operators, consistently with a previous introductory paper by the authors in \cite{kapidaniHighOrderGeometric2022}. We also present the definition of the primal and the dual spline complexes, which will be used in the discretization of the problem.

\subsection{Maxwell's equations and differential forms}
From here onwards let us set the space dimension $n=3$ and the domain $\Omega\subset\mathbb{R}^3$. For non-negative integers $k$, we denote the space of smooth differential $k$--forms with $\Lambda^k(\Omega)$. In general, considering smooth functions is too restrictive. We will instead need the Hilbert space $L^2\Lambda^k(\Omega)$, defined as the completion of $\Lambda^k(\Omega)$ with respect to the $L^2$--inner product, see \cite{AFW06}.

The exterior derivative $d^k$ maps $k$--forms into $(k+1)$--forms. 
From here onwards we will simply denote the exterior derivative by $d$ (as done in \eqref{eq:ampere}--\eqref{eq:divgaugeB}) when there is no confusion on the order of the differential form. An important property of the exterior derivative is that $d \circ d \omega = 0$ for any differential $k$--form $\omega$.
Following \cite{AFW06}, let us define the Sobolev spaces, for $k= 0, \ldots, 3$,
\begin{equation*}
H\Lambda^k(\Omega) = \left\{ \omega \in L^2\Lambda^k(\Omega) : d^k \omega \in L^2\Lambda^{k+1}(\Omega) \right\},
\end{equation*}
from which we can construct the $L^2$ de Rham complex of differential forms
\begin{equation} \label{eq:derhamcomplex}
\begin{CD}
H\Lambda^0(\Omega) @>d^0>> H\Lambda^1(\Omega) @>d^1>> H\Lambda^2(\Omega) @>d^2>> H\Lambda^3(\Omega).
\end{CD}
\end{equation}
Furthermore, we will also need to define de Rham complexes of differential forms with vanishing boundary traces, which are built from spaces of differential $k$--forms with compact support. We distinguish them from the ones in \eqref{eq:derhamcomplex} by using the zero subscript, and we build the sequence
\begin{equation} \label{eq:derhamcomplex_homogeneous}
\begin{CD}
H_0\Lambda^0(\Omega) @>d^0>> H_0\Lambda^1(\Omega) @>d^1>> H_0\Lambda^2(\Omega) @>d^2>> H_0\Lambda^3(\Omega).
\end{CD}
\end{equation}

Assuming for simplicity that the domain $\Omega$ is bounded by a perfect electrical conductor, the tangential component of the electric field vanishes on $\partial \Omega$, and therefore we have to solve \eqref{eq:ampere}--\eqref{eq:divgaugeB} for differential 1--forms $\mathrm{E} \in H_0 \Lambda^1(\Omega)$ and $\mathrm{H} \in H \Lambda^1(\Omega)$, and differential 2--forms $\mathrm{B} \in H_0 \Lambda^2(\Omega)$, $\mathrm{D} \in H \Lambda^2(\Omega)$. Thus, we will approximate the differential forms $\mathrm{E,B}$ by discrete differential forms that belong to subspaces of the sequence \eqref{eq:derhamcomplex_homogeneous}, and $\mathrm{H,D}$ using discrete subspaces of the sequence \eqref{eq:derhamcomplex}. The system of equations is completed with initial conditions $\mathrm{D}(\boldsymbol{x},0) = \mathrm{D}_0 (\boldsymbol{x})$, and $\mathrm{B}(\boldsymbol{x},0) = \mathrm{B}_0 (\boldsymbol{x})$.

\subsection{Hodge star operators}
To perform discrete time integration on the discrete Maxwell system, the unknowns on the two sequences must be complemented by constitutive equations, which should map 1--forms into 2--forms, or vice versa, and which take into account the material properties. This is expressed in terms of Hodge star operators, as in \eqref{eq:cont_hodge1}, that we now define precisely.

Let us assume that $\gamma$ is a bounded and uniformly positive scalar-valued field defined in $\Omega$. The Hodge star operator $\star_{\gamma}^k:L^2 \Lambda^k(\Omega) \rightarrow L^2 \Lambda^{3-k}(\Omega)$ is a linear operator from $k$--forms into $(3-k)$--forms.
One property is sufficient to completely define it for any $\omega\in L^2\Lambda^k(\Omega)$, namely the following equality:
\begin{equation} \label{eq:def_hodge_diff}
    \left( \eta, \omega \right)_{L_\gamma^2 \Lambda^k(\Omega)} = \int_{\Omega} \eta \wedge \star_{\gamma}^k(\omega),
    \quad \forall \eta \in L^2 \Lambda^k(\Omega), 
\end{equation}
where $(\cdot, \cdot)_{L^2_\gamma\Lambda^k}$ is the inner product in $L^2$ with the material parameter $\gamma$ as its weight, and where we have introduced the exterior product (or wedge product, denoted by $\wedge$), an alternating bilinear product between a differential $k$--form and a differential $l$--form, yielding a differential $(k+l)$--form. By virtue of this, we always obtain a 3--form integrated on the 3--manifold $\Omega$ on the right-hand side of \eqref{eq:def_hodge_diff}. An important property of the Hodge star operator is its invertibility, and in fact for dimension $n=3$ it holds that $\star^k_\gamma \circ \star^{3-k}_{1/\gamma} = \mathrm{Id}$. 
Alternatively, and thanks to the invertibility of the operator, the constitutive equations can be also written as
\begin{equation} \label{eq:cont_hodge2}
    \begin{aligned}
    & \star^2_{1/\mu} \mathrm{B} = \mathrm{H},
    & \star^2_{1/\epsilon} \mathrm{D} = \mathrm{E}.
    \end{aligned}
    \end{equation}

Finally, by combining \eqref{eq:ampere}, \eqref{eq:faraday} and \eqref{eq:cont_hodge1} with the property \eqref{eq:def_hodge_diff}, it is straightforward to derive the exterior calculus notation for the conserved energy of the Maxwell system as 
\begin{equation}\label{eq:em_energy_cont}
    \mathcal{E} = \frac{1}{2} \left(\int_\Omega \mathrm{E} \wedge \mathrm{D} + \int_\Omega \mathrm{H} \wedge \mathrm{B} \right).
\end{equation}

\subsection{Spline complexes of differential forms} \label{sec:splines}

We will work with spline discretizations as standard in IGA, and we will make the assumption that the domain is {described} by a single patch, i.e., we define the physical domain $\Omega \subset \mathbb{R}^3$ through a parametrization of the form $\textbf{F}: \hat \Omega \rightarrow \Omega$, where $\hat \Omega = (0,1)^3$ is called the parametric domain. We thus relate differential $k$--forms in the parametric domain to differential $k$--forms in the physical domain using a set of pullback operators $\iota^k: H \Lambda^k(\Omega) \rightarrow H\Lambda^k(\hat \Omega)$, their expression in the three-dimensional case can be found for instance in \cite{BRSV11}, in terms of vector proxies. 
An important property is that the pullback commutes both with the exterior derivative and the wedge product.
These tools allow us to define the primal and dual complex of splines, as introduced in \cite{Buffa_2020aa}. We refer to it for more details, and to \cite{kapidaniHighOrderGeometric2022} for a presentation in terms of differential forms. For the purposes of the present work we again restrict ourselves to the three-dimensional case.

Let $p>1$ denote the polynomial degree of univariate B-splines, we introduce the $p$-open knot vector $\Xi = \{\xi_{1}, \ldots , \xi_{m+p+1} \}$, where $m$ is the number of basis functions. We denote by $S_p(\Xi)$ the space spanned by them, which is the space of piecewise polynomials of degree $p$ with $p-r_i$ continuous derivatives at each knot $\xi_i$, where $r_i$ is the multiplicity of the knot. Assuming that the multiplicity is $r_i < p-1$ for every internal knot, the functions in $S_p(\Xi)$ are at least $C^1$ continuous. Their derivatives then belong to the space $S_{p-1}(\Xi')$, with $\Xi' = \{\xi_2, \ldots, \xi_{m+p}\}$ defined from $\Xi$ by removing the first and last repeated knots. Analogously, functions in $S_{p-1}(\Xi')$ are at least $C^0$ continuous, and their (weak) derivatives belong to the space $S_{p-2}(\Xi'')$, with $\Xi'' = \{\xi_3, \ldots, \xi_{m+p-1}\}$.

Multivariate B-splines are defined by tensor product. Let us assume for simplicity that the same degree $p$ and knot vector $\Xi$ are used in the three directions of the space. Then, one can construct a discrete de Rham complex of B-splines of the form
\begin{equation*} 
\begin{CD}
\Xhath{0} @>d^0>> \Xhath{1} @>d^1>> \Xhath{2} @>d^{2}>> \Xhath{3},
\end{CD}
\end{equation*}
where the spaces are defined by suitable tensor-products of univariate spaces $S_p(\Xi)$ and $S_{p-1}(\Xi')$. 
A discrete complex for spaces with homogeneous boundary conditions, that we denote by $\Xhathb{k}$, is obtained analogously, removing the first and last basis functions of the univariate space $S_p(\Xi)$, and leaving $S_{p-1}(\Xi')$ unchanged.

To define the spaces of the dual complex we proceed in a completely analogous fashion, replacing splines of degrees $p$ and $p-1$ by splines of degree $p-1$ and $p-2$, respectively, to build a dual complex on the parametric domain, whose elements will be accordingly labelled by $\XhathT{k}$. Finally, obtaining the discrete sequences in the physical domain is a matter of applying the correct pullback operators, namely
\begin{equation*}
X_{h,0}^k := \left\{ \omega : \iota^k(\omega) \in \Xhathb{k} \right\}, \quad
\tilde{X}_h^k := \left\{ \omega : \iota^k(\omega) \in \XhathT{k} \right\}.
\end{equation*}
To ensure the expected approximation properties for all the spaces in the physical domain, it is necessary for the parameterization $\textbf{F}$ to have at least the same regularity as the corresponding spaces in the parametric domain. Since the space with highest regularity is $\hat{X}_h^0$, the requirement can be easily met by defining each component of the map $\mathbf{F}$ as a (rational) spline living in said space. We refer to \cite{BRSV11} for more technical details.

With these definitions, we have obtained a primal complex of splines of mixed degree $p$ and $p-1$, for spaces with homogeneous boundary conditions $X^k_{h,0} \subset H_0\Lambda^k(\Omega)$,
\begin{equation} \label{eq:cd-iga-tp-drchlt}
    \begin{CD}
    {X}^{0}_{h,0} @>d^0>> {X}^{1}_{h,0} @>d^1>> {X}^{2}_{h,0} @>d^{2}>> {X}^{3}_h,
    \end{CD}
\end{equation}
and a dual complex of spline spaces of mixed degree $p-1$ and $p-2$, given by 
\begin{equation} \label{eq:cd-iga-tp-nmnn}
    \begin{CD}
    \XhT{0} @>d^0>> \XhT{1} @>d^1>> \XhT{2} @>d^{2}>> \XhT{3},
    \end{CD}
\end{equation}
where this time $\XhT{k} \subset H\Lambda^k(\Omega)$, and it has been proven in \cite{hiemstra_isogeometric_2011,Buffa_2020aa} that $\dim X^k_{h,0} = \dim \XhT{3-k}$.

In light of the above definitions, and recalling that we assume boundary conditions for a perfect electrical conductor, we discretize the electric field as a 1--form of the primal complex, $\mathrm{E}_h \in X^1_{h,0}$, the magnetic field as a 1--form of the dual complex, $\mathrm{H}_h \in \XhT{1}$, \textcolor{black}{and correspondingly the electric displacement is discretized as a 2--form of the dual complex, $\mathrm{D}_h \in \XhT{2}$, and the magnetic induction as a 2--form of the primal complex, $\mathrm{B}_h \in X^2_{h,0}$.}

\subsection{Pairing matrices}
In the discrete setting, given bases for the discrete spaces, we can define pairing matrices between the spaces of primal differential $k$--forms $X^k_{h,0}$ and dual differential $(3-k)$--forms $\XhT{3-k}$. 
Let $\omega_h \in X^k_{h,0}$ and $\tilde{\eta}_h \in \XhT{3-k}$ be respectively represented by the vectors of degrees of freedom $\boldsymbol{\omega}$ and $\tilde{\boldsymbol{\eta}}$. The square pairing matrices ${\bf K}_k$ and $\tilde{\bf K}_{3-k}$, are respectively determined by
\begin{align*}
\tilde{\boldsymbol{\eta}}^\top {\bf K}_k \boldsymbol{\omega} = \int_\Omega \tilde{\eta}_h \wedge \omega_h,  \quad \text{ for all } \omega_h \in X^k_{h,0}, \, \tilde{\eta}_h \in \XhT{3-k} , \\
\boldsymbol{\omega}^\top \tilde{\bf K}_{3-k} \tilde{\boldsymbol{\eta}} = \int_\Omega \omega_h \wedge \tilde{\eta}_h, \quad \text{ for all } \omega_h \in X^k_{h,0}, \, \tilde{\eta}_h \in \XhT{3-k}, 
\end{align*}
i.e., the matrix entries are integrals of wedge products between basis elements of the two spaces of forms.
From the properties of the wedge product, it immediately follows in the three-dimensional case that 
\begin{equation} \label{eq:pairing_transpose}
    \tilde{\bf K}_{3-k} = (-1)^{k (3-k)} ({\bf K}_k)^\top = ({\bf K}_k)^\top.
\end{equation}
Since the pairing $(X^k_{h,0}, \XhT{3-k})$ is stable, as proved in the framework of mortar methods in \cite{Buffa_2020aa}, all the pairing matrices are invertible.  The properties of the wedge product {guarantee that} the pairing matrices are metric independent, and thus their inverse, or their LU factorization, can be computed in the parametric domain $\hat \Omega$, exploiting the tensor-product structure of splines. We will see more details in Section~\ref{sec:tensorized}.

To conclude the section, we remark that with a suitable choice of the basis functions of univariate spline spaces, the exterior derivative of discrete splines can be written in terms of incidence matrices associated to a Cartesian mesh, see \cite{ratnani} and \cite{BSV14} for more details. 
We will denote these matrices for the primal and the dual complex respectively by $\mathbf{D}^k$ and $\tilde{\mathbf{D}}^k$, for $k=0,1,2$.
There is a close relation between the pairing matrices and the incidence matrices of the exterior derivative. Indeed, from their corresponding definitions, and recalling that the functions in the primal complex have vanishing boundary conditions, for any $k>0$ it holds that 
\begin{equation} \label{eq:disc_int_by_parts1}
    ({\bf D}^{k-1})^\top \tilde{\bf K}_{3-k} = (-1)^{k} \tilde{\bf K}_{3-k+1} \tilde{\bf D}^{3-k},
\end{equation}
and conversely, using \eqref{eq:pairing_transpose} it holds that
\begin{equation} \label{eq:disc_int_by_parts2}
    (\tilde{\bf D}^{3-k})^\top {\bf K}_{k-1} = (-1)^{3-k+1} {\bf K}_{k} {\bf D}^{k-1},
\end{equation}
which both represent the discrete counterpart of the integration by parts formula
\begin{equation*} 
    \int_\Omega d\omega \wedge \eta = 
(-1)^k \int_\Omega \omega \wedge d \eta, \quad \text{ for } \omega \in H_0 \Lambda^{k-1}(\Omega), \eta \in H\Lambda^{3-k} (\Omega), 
\end{equation*}
\noindent for $0 < k \leq 3$.
It is important to remark that the metric does not play any role in either definition of pairing matrices or incidence matrices, and consequently also in none of the properties above.

\section{Spatial discretization with energy preserving Hodge--star operators}\label{sec:hodge}
In this section we present the discretization in space of Maxwell's equations using splines of the primal and dual complex above, along with two possible choices for the discrete Hodge star operators which yield two different discretization schemes. {The first scheme we present involves the solution of linear systems for mass matrices, and is equivalent to standard Galerkin techniques. The second scheme is the main contribution of this paper: it requires the solution of linear systems for pairing matrices, and we prove that it is equivalent to a Petrov-Galerkin scheme.} We show that both choices of the Hodge operators provide a semi-discrete method which preserves electromagnetic energy across time, and which also satisfies the two Gauss laws for the conservation of charges.
From now on we will assume that in Maxwell's equations all the material parameters are time-invariant, which implies that Hodge operators commute with time derivatives. 

\subsection{Spatial discretization of Amp\`ere-Maxwell and Faraday equations}
We start recalling that the four unknowns of the semi-discrete system are $\mathrm{E}_h \in X^1_{h,0}$, $\mathrm{H}_h \in \XhT{1}$, $\mathrm{D}_h \in \XhT{2}$, $\mathrm{B}_h \in X^2_{h,0}$, and we will denote their corresponding vectors of degrees of freedom by  $\mathbf{e}$, $\mathbf{h}$, $\mathbf{d}$ and $\mathbf{b}$, respectively. We will also assume that the known source is a discrete field $\mathrm{J}_h \in \tilde{X}^2_h$ which satisfies $d \mathrm{J}_h = 0$, and denote the corresponding vector of coefficients as $\mathbf{j}$. With the choice made above for the discrete differential forms, and recalling that the exterior derivative for spline spaces can be written in terms of the incidence matrices $\mathbf{D}^k$ and $\tilde{\mathbf{D}}^k$, the discrete version of Maxwell's equations \eqref{eq:ampere}--\eqref{eq:faraday} is written as
\begin{equation} \label{eq:maxwell_discrete}
\begin{aligned}
& \partial_t \mathbf{d} = \tilde{\mathbf{D}}^1 \mathbf{h} - \mathbf{j}, \\
& \partial_t \mathbf{b} = - \mathbf{D}^1 \mathbf{e}.
\end{aligned}
\end{equation}
It is important to remark that only the first two equations in \eqref{eq:ampere}--\eqref{eq:divgaugeB} need to be discretized. Indeed, the last two equations are automatically satisfied by the fact that applying the exterior derivative twice always vanishes, also in the discrete case, and since the two Gauss laws are satisfied exactly, this implies the exact conservation of charges.

The electromagnetic energy \eqref{eq:em_energy_cont} is approximated in the discrete setting by the quantity 
\begin{equation*}
\mathcal{E}_h = \frac{1}{2} \left(\int_\Omega \mathrm{E}_h \wedge \mathrm{D}_h + \int_\Omega \mathrm{H}_h \wedge \mathrm{B}_h \right),
\end{equation*}
that we can represent in matrix form with the help of the pairing matrices as
\begin{equation} \label{eq:energy_matrix_wedge}
\mathcal{E}_h = \frac{1}{2} \left( \mathbf{e}^\top \tilde{\mathbf{K}}_2 \mathbf{d} + \mathbf{h}^\top \mathbf{K}_2 \mathbf{b} \right).
\end{equation}
To complete the spatial discretization, we are only left with the need to approximate the Hodge star operators.

\subsection{Discrete Hodge star operators}
We now present two different alternatives for the discrete Hodge star operators, 
which particularize to our setting the general construction in \cite{HIP99c}. The two choices, which depend on whether we discretize \eqref{eq:cont_hodge1} or \eqref{eq:cont_hodge2}, will lead to two different numerical schemes. Let us introduce first the discrete version of the Hodge star operators in \eqref{eq:cont_hodge1}, for which we have to define two operators of the form
\[
\star^1_{h,\epsilon} : X_{h,0}^1 \longrightarrow \tilde{X}_h^2,
\qquad 
\tilde{\star}^1_{h,\mu} : \tilde{X}_{h}^1 \longrightarrow X_{h,0}^2,
\]
i.e., one operator from primal 1--forms into dual 2--forms and one operator from dual 1--forms into primal 2--forms, where we use the $\tilde{\cdot}$ notation in the discrete Hodge to highlight that its domain of definition is a discrete space of the dual sequence. 

Let $\mathrm{E}_h \in X^1_{h,0}$, then the application of the discrete Hodge star operator $\star_{{h,\epsilon}}^1$ on $\mathrm{E}_h$ mimics the definition of the continuous Hodge star operator in \eqref{eq:def_hodge_diff}. By setting $\mathrm{D}_h = \star_{{h,\epsilon}}^1 \mathrm{E}_h$, this is uniquely determined by
\begin{equation} \label{eq:def_hodge_discrete1}
    \int_\Omega \eta_h \wedge \mathrm{D}_h = \left( \eta_h, \mathrm{E}_h \right)_{L_\epsilon^2 \Lambda^1(\Omega)} ,
    \quad \forall \eta_h \in X_{h,0}^1,
\end{equation}
and the definition of the second Hodge star operator $\tilde{\star}^1_{h,\mu}$, for $\mathrm{B}_h = \tilde{\star}^1_{h,\mu} \mathrm{H}_h$,  is completely analogous, replacing $X^1_{h,0}$ by $\tilde{X}^1_h$.

Alternatively, we can introduce a discrete version of the Hodge star operators \eqref{eq:cont_hodge2}, which means that we will define two operators of the form
\[
\tilde{\star}^2_{h,1/\epsilon} : \tilde{X}_h^2 \longrightarrow X_{h,0}^1,
\qquad 
\star^2_{h,1/\mu} : X_{h,0}^2 \longrightarrow \tilde{X}_{h}^1,
\]
i.e., the operators now map 2--forms into 1--forms. By mimicking again the definition of the continuous Hodge star operator, given $\mathrm{D}_h \in \XhT{2}$, and setting $\mathrm{E}_h = \tilde{\star}^2_{h,1/\epsilon} \mathrm{D}_h$, the discrete Hodge star operator $\tilde{\star}^2_{h,1/\epsilon}$ is uniquely determined by
\begin{equation}\label{eq:def_hodge_discrete2}
    \int_\Omega \eta_h \wedge \mathrm{E}_h = \left( \eta_h, \mathrm{D}_h \right)_{L_{1/\epsilon}^2 \Lambda^2(\Omega)} ,
    \quad \forall \eta_h \in \XhT{2},
\end{equation}
and similarly the second discrete operator $\star^2_{h,1/\mu}$, for $\mathrm{H}_h = \star^2_{h,1/\mu} \mathrm{B}_h$, is defined through replacing $\XhT{2}$ by $X_{h,0}^2$.


\begin{remark}
The two kinds of operators defined here are inspired by the work of Hiptmair \cite{HIP99c}, and correspond to the two global operators from the primal to the dual complex analysed in \cite{kapidaniHighOrderGeometric2022}. Unfortunately, the third, local operator studied in the same paper does not seem to maintain the same approximation properties when mapping from the dual to the primal complex, and would not guarantee conservation of electromagnetic energy. We have therefore decided to discard it in the present paper.
\end{remark}

\subsection{First discretization scheme} \label{sec:first_hodge}


The first scheme considers the discrete operators $\star^1_{h,\epsilon}$ and $\tilde \star^1_{h,\mu}$. Using \eqref{eq:def_hodge_discrete1} in its matrix form, with boldface symbols for vectors of degrees of freedom, the definition of the operator $\star^1_{h,\epsilon}$ gives
\[
{\boldsymbol{\eta}^\top \tilde{\bf K}_{2} \mathbf{d} = \boldsymbol{\eta}^\top  {\bf M}^1_\epsilon \mathbf{e},}
\]
where ${\bf M}^1_\epsilon$ is the standard mass matrix for $X^1_{h,0}$ and $\tilde{\bf K}_{2}$ is the pairing matrix between $X^1_{h,0}$ and $\XhT{2}$. Since the vector $\boldsymbol{\eta}$ is arbitrary, the equation can be simply written as
\begin{equation*}
{\bf M}^1_\epsilon \mathbf{e} = \tilde{\bf K}_{2} \mathbf{d}.
\end{equation*}

Analogously, by setting $\mathrm{B}_h = \tilde{\star}^1_{h,\mu} \mathrm{H}_h$, following the same reasoning, the role of the discrete Hodge operator $\tilde{\star}^1_{h,\mu}$ is expressed in matrix form as
\begin{equation*}
\tilde{{\bf M}}^1_\mu \mathbf{h} = \mathbf{K}_{2} \mathbf{b},
\end{equation*}
where $\tilde{\mathbf{M}}^1_\mu$ is the mass matrix corresponding to $\tilde{X}^1_h$.
Finally, combining \eqref{eq:maxwell_discrete} with the definition of the two discrete Hodge star operators, the semi-discrete scheme after spatial discretization is given by
\begin{align}
& \partial_t \mathbf{d} = \tilde{\mathbf{D}}^1 \mathbf{h} - \mathbf{j}, \label{eq:discrete1_ampere}\\
& {\bf M}^1_\epsilon \mathbf{e} = \tilde{\bf K}_{2} \mathbf{d}, \label{eq:discrete1_hodge_eps}\\
& \partial_t \mathbf{b} = - \mathbf{D}^1 \mathbf{e}, \label{eq:discrete1_faraday}\\
& \tilde{{\bf M}}^1_\mu \mathbf{h} = \mathbf{K}_{2} \mathbf{b}, \label{eq:discrete1_hodge_mu}
\end{align}
where, before dealing with time integration, the first and third equations only require the application of the incidence matrices. The second and fourth equations, which represent the application of the discrete Hodge star operators, instead require a matrix-vector multiplication with a pairing matrix, and the solution of a linear system for a mass matrix.

\subsubsection{Energy conservation}
For the study of conservation of energy, we assume that the source current $\mathrm{J}$ is equal to zero. Then, to prove that the discrete energy is conserved by the scheme \eqref{eq:discrete1_ampere}--\eqref{eq:discrete1_hodge_mu}, let us sum up equation \eqref{eq:discrete1_ampere} multiplied by $\mathbf{e}^\top \tilde{\mathbf{K}}_2$ and equation \eqref{eq:discrete1_faraday} multiplied by $\mathbf{h}^\top \mathbf{K}_2$, to obtain
\[
\mathbf{e}^\top \tilde{\mathbf{K}}_2 \partial_t \mathbf{d} + \mathbf{h}^\top \mathbf{K}_2 \partial_t \mathbf{b} = \mathbf{e}^\top \tilde{\mathbf{K}}_2 \tilde{\mathbf{D}}^1 \mathbf{h} - \mathbf{h}^\top \mathbf{K}_2 \mathbf{D}^1 \mathbf{e}.
\]
Then, by applying the property of integration by parts in \eqref{eq:disc_int_by_parts2} with $k=2$, and the property of the pairing matrices \eqref{eq:pairing_transpose}, we get
\begin{equation} \label{eq:energy_aux}
\begin{aligned}
\mathbf{e}^\top \tilde{\mathbf{K}}_2 \partial_t \mathbf{d} + \mathbf{h}^\top \mathbf{K}_2 \partial_t \mathbf{b} 
& = \mathbf{e}^\top \tilde{\mathbf{K}}_2 \tilde{\mathbf{D}}^1 \mathbf{h} - {\mathbf{h}^\top (\tilde{\mathbf{D}}^1)^\top} \mathbf{K}_1 \mathbf{e} \\
& = 
\mathbf{e}^\top \tilde{\mathbf{K}}_2 \tilde{\mathbf{D}}^1 \mathbf{h} - {\mathbf{h}^\top (\tilde{\mathbf{D}}^1)^\top \tilde{\mathbf{K}}_2^\top} \mathbf{e} = 0.
\end{aligned}
\end{equation}
At this point we use the definition of the discrete Hodge star operators. Recalling that the pairing matrices $\mathbf{K}_2$ and $\tilde{\mathbf{K}}_2$ are invertible, we can replace $\mathbf{d}$ and $\mathbf{b}$ by their respective expressions from \eqref{eq:discrete1_hodge_eps} and \eqref{eq:discrete1_hodge_mu}, from which we get
\begin{equation*}
\mathbf{e}^\top \partial_t (\mathbf{M}^1_\epsilon \mathbf{e}) + \mathbf{h}^\top \partial_t \tilde{\mathbf{M}}^1_{\mu} \mathbf{h}  = 0.
\end{equation*}
Since we have assumed that both $\epsilon$ and $\mu$ are time-invariant, the mass matrices can be factored out of the partial derivative, and the conservation of the energy follows by standard arguments.

It is important to note that we can also replace $\mathbf{d}$ and $\mathbf{b}$ in \eqref{eq:energy_matrix_wedge} by their expressions in terms of the discrete Hodge star operators, to get an equivalent expression for the energy:
\begin{equation*}
\mathcal{E}_h = \frac{1}{2} \left( \mathbf{e}^\top \mathbf{M}^1_\epsilon \mathbf{e} + \mathbf{h}^\top \tilde{\mathbf{M}}^1_\mu \mathbf{h} \right),
\end{equation*}
i.e., the standard quadratic form based on mass matrices considered in finite element schemes which require the solution of linear systems for the same symmetric positive definite mass matrices.


\subsubsection{Weak formulation in terms of $\mathrm{E}$ and $\mathrm{H}$} \label{sec:weak_first}
We have presented the method using four unknowns and four equations, where equations involving the time derivatives and exterior derivatives are solved in their strong form, while the Hodge star operators are imposed weakly. This is akin to the usual procedure in FIT formulations, see e.g. the formulation in \cite{codecasa_explicit_2008}. Alternatively, it is possible to write a weak formulation of the problem only in terms of $\mathrm{E}$ and $\mathrm{H}$, as it is usually done in finite elements. We start by multiplying equations \eqref{eq:discrete1_ampere} and \eqref{eq:discrete1_faraday} by the same terms used to prove energy conservation, and replacing the expressions of $\mathbf{d}$ and $\mathbf{b}$ via application of Hodge star operators. Since we assume that material properties are independent of time, we then obtain
\begin{align*}
& {\bf M}^1_\epsilon \partial_t \mathbf{e} = \tilde{\bf K}_{2} \tilde{\mathbf{D}}^1 \mathbf{h} - \tilde{\bf K}_{2} \mathbf{j}, \\
& \tilde{{\bf M}}^1_\mu \partial_t \mathbf{h} = - \mathbf{K}_{2} \mathbf{D}^1 \mathbf{e}.
\end{align*}

By the definition of mass and pairing matrices, and using the fact that incidence matrices encode the action of exterior derivatives, the above is equivalent to the weak formulation
\begin{equation*}
\begin{aligned}
& (\partial_t \mathrm{E}_h, \eta_h)_{L^2_\epsilon\Lambda^1(\Omega)} = \int_\Omega \eta_h \wedge (d \mathrm{H}_h - \mathrm{J}_h) \; & \text{ for all }\eta_h \in X^1_{h,0},\\
& (\partial_t \mathrm{H}_h, \xi_h)_{L^2_\mu\Lambda^1(\Omega)} = - \int_\Omega \xi_h \wedge d \mathrm{E}_h \; & \text{ for all }\xi_h \in \tilde{X}^1_{h}.
\end{aligned}
\end{equation*}
Thus, our first discretization scheme is equivalently written as a Galerkin method with a mixed formulation, {similar for instance to \cite{ratnani}}, with the difference that $\mathrm{E}_h$ is discretized in the primal complex and $\mathrm{H}_h$ is discretized in the dual complex.

\subsection{Second discretization scheme} \label{sec:solve_K}

Alternatively with respect to the previous subsection, we can {obtain a different discretization scheme by considering} the Hodge operators $\tilde \star^2_{h,1/\epsilon}$ and $\star^2_{h,1/\mu}$ mapping 2--forms into 1--forms, such as the one defined by \eqref{eq:def_hodge_discrete2}. 
Similar arguments as above show that the first operator can be equivalently written in matrix form as
\begin{equation*}
\mathbf{K}_1 \mathbf{e} = \tilde{\mathbf{M}}^2_{1/\epsilon} \mathbf{d},
\end{equation*}
where $\tilde{\mathbf{M}}^2_{1/\epsilon}$ is the mass matrix of the space of dual 2--forms $\XhT{2}$, and $\mathbf{K}_1$ is a pairing matrix between $\XhT{2}$ and $X^1_{h,0}$. In an analogous way, given $\mathrm{B}_h \in X^1_{h,0}$ we compute $\mathrm{H}_h = \star^2_{h,1/\mu} \mathrm{B}_h$ in terms of their degrees of freedom as 
\begin{equation*}
\tilde{\mathbf{K}}_1 \mathbf{h} = \mathbf{M}^2_{1/\mu} \mathbf{b},
\end{equation*}
where $\mathbf{M}^2_{1/\mu}$ is the mass matrix of the space of 2--forms $X^2_{h,0}$.
As we did for the first scheme, we combine the definition of the discrete Hodge star operators with \eqref{eq:maxwell_discrete}, to obtain the spatial semidiscretization 
\begin{align}
& \partial_t \mathbf{d} = \tilde{\mathbf{D}}^1 \mathbf{h} - \mathbf{j}, \label{eq:discrete2_ampere}\\
& \mathbf{K}_1 \mathbf{e} = \tilde{\mathbf{M}}^2_{1/\epsilon} \mathbf{d}, \label{eq:discrete2_hodge_eps}\\
& \partial_t \mathbf{b} = - \mathbf{D}^1 \mathbf{e}, \label{eq:discrete2_faraday}\\
& \tilde{\mathbf{K}}_1 \mathbf{h} = \mathbf{M}^2_{1/\mu} \mathbf{b}, \label{eq:discrete2_hodge_mu}
\end{align}
\textcolor{black}{where the equations look very similar to the ones for the first scheme in \eqref{eq:discrete1_ampere}--\eqref{eq:discrete1_hodge_mu}. There is nevertheless a very important difference due to the choice of the Hodge star operators: the second and fourth equations require the solution of linear systems associated to the pairing matrices, instead of mass matrices. At first glance this might come across as a drawback, since these matrices are {in general} not symmetric positive definite. There are however important advantages, coming from the tensor-product structure of B-splines and the fact that the pairing matrices are metric-independent, as we will see in Section~\ref{sec:tensorized}.}

\subsubsection{Energy conservation}
Regarding energy conservation, we remark that equation \eqref{eq:energy_aux} was obtained without making use of discrete Hodge star operators. It is thus still valid in this case, yielding 
\begin{align*}
\mathbf{e}^\top \tilde{\mathbf{K}}_2 \partial_t \mathbf{d} + \mathbf{h}^\top \mathbf{K}_2 \partial_t \mathbf{b} = \mathbf{e}^\top \tilde{\mathbf{K}}_2 \tilde{\mathbf{D}}^1 \mathbf{h} - \mathbf{h}^\top \mathbf{K}_2 \mathbf{D}^1 \mathbf{e} = 0.
\end{align*}

We now make use of the discrete Hodge star operators and replace $\mathbf{e}$ and $\mathbf{h}$ by their respective expressions in \eqref{eq:discrete2_hodge_eps} and \eqref{eq:discrete2_hodge_mu}. Using the symmetry of mass matrices and the property of the pairing matrices~\eqref{eq:pairing_transpose}, we obtain
\begin{align*}
{((\mathbf{K}_1)^{-1} \tilde{\mathbf{M}}^2_{1/\epsilon} \mathbf{d})^\top \tilde{\mathbf{K}}_2 \partial_t \mathbf{d} + (\tilde{\mathbf{K}}_1^{-1} \mathbf{M}^2_{1/\mu} \mathbf{b})^\top \mathbf{K}_2 \partial_t \mathbf{b}} =
\mathbf{d}^\top \tilde{\mathbf{M}}^2_{1/\epsilon} \partial_t \mathbf{d} + \mathbf{b}^\top \mathbf{M}^2_{1/\mu} \partial_t \mathbf{b} = 0,
\end{align*}
and the conservation of energy follows from the fact that the mass matrices are symmetric positive definite.

Analogously to what we saw for the first scheme, replacing the discrete Hodge star operators in \eqref{eq:energy_matrix_wedge} gives an equivalent expression of the energy in terms of mass matrices, which is given by
\begin{equation*}
\mathcal{E}_h = \mathbf{d}^\top \tilde{\mathbf{M}}^2_{1/\epsilon} \mathbf{d} + \mathbf{b}^\top \mathbf{M}^2_{1/\mu} \mathbf{b},
\end{equation*}
which again yields a quadratic form based on mass matrices, in this case for the spaces of discrete 2--forms.

\subsubsection{Weak formulation in terms of $\mathrm{E}$ and $\mathrm{H}$}
Analogously to the first scheme, it is possible to write the discrete problem in an equivalent weak formulation. Multiplying equations \eqref{eq:discrete2_ampere} and \eqref{eq:discrete2_faraday} respectively by $\tilde{\mathbf{M}}^2_{1/\epsilon}$ and $\tilde{\mathbf{M}}^2_{1/\mu}$, and replacing $\mathbf{d}$ and $\mathbf{b}$ by their respective expressions in the definition of the Hodge star operators, we get
\begin{align*}
& {\bf K}_1 \partial_t \mathbf{e} = \tilde{\mathbf{M}}^2_{1/\epsilon} \tilde{\mathbf{D}}^1 \mathbf{h} - \tilde{\mathbf{M}}^2_{1/\epsilon} \mathbf{j}, \\
& \tilde{\bf K}_1 \partial_t \mathbf{h} = - \mathbf{M}^2_{1/\mu} \mathbf{D}^1 \mathbf{h},
\end{align*}
which, through the same arguments as above, is equivalent to the weak formulation
\begin{equation*}
\begin{aligned}
& \int_\Omega \eta_h \wedge \partial_t \mathrm{E}_h = (d \mathrm{H}_h - \mathrm{J}_h, \eta_h)_{L^2_{1/\epsilon} \Lambda^2(\Omega)}  \; & \text{ for all }\eta_h \in \tilde{X}^2_{h}, \\
& \int_\Omega \xi_h \wedge \partial_t \mathrm{H}_h = - (d \mathrm{E}_h, \xi_h)_{L^2_{1/\mu} \Lambda^2(\Omega)} \; & \text{ for all }\xi_h \in X^2_{h,0},
\end{aligned}
\end{equation*}
through which one can conclude that our second discretization scheme is equivalent to a Petrov-Galerkin formulation of the continuous problem, where the unknowns $\mathrm{E}_h$ and $\mathrm{H}_h$ are 1--forms, respectively defined in the primal and the dual complex. Their corresponding spaces for test functions are the spaces of 2--forms from the other complex.


\begin{remark}
In most of the methods of the literature based on a dual grid \cite{Clemens-Weiland,hirani}, applying the Hodge operator is identified with solving a linear system for the (lumped) mass matrix. As pointed out by Hiptmair \cite{HIP99c}, this is true in the particular case in which the pairing matrices coincide exactly with the identity matrix, a condition which can be fulfilled only for particular choices of the dual mesh. As far as we know, methods based on a diagonal system matrix have never been extended to high order basis functions for the Maxwell equations. 
\end{remark}

\subsection{Generalization to non-homogeneous boundary conditions}
As we have seen, the two methods can be understood as Galerkin or Petrov-Galerkin schemes. Therefore, to impose a non-homogeneous boundary condition for the electric field, the procedure is the same as for a standard Galerkin technique. We write the electric field as a linear combination of two components $\mathrm{E}_h = \mathrm{E}_{h,0} + \mathrm{E}_{h,b}$, where $\mathrm{E}_{h,0} \in X^1_{h,0}$ has the {tangential} trace vanishing on $\partial\Omega$ and $\mathrm{E}_{h,b} \in X^1_h$ is a suitable lifting, which can be computed through a local projection of the boundary conditions on its trace space. This computation only involves the basis functions of $X^1_h$ with a non-vanishing boundary trace, see for instance \cite{Vazquez_2016aa}.

Once the lifting $\mathrm{E}_{h,b}$ has been computed, it enters the equations into the right-hand side. In particular, replacing the expression of $\mathrm{E}_h$ in \eqref{eq:def_hodge_discrete1} and rearranging terms, we obtain
\begin{equation*}
    \left( \eta_h, \mathrm{E}_{h,0} \right)_{L_\epsilon^2 \Lambda^1(\Omega)} = \int_\Omega \eta_h \wedge \mathrm{D}_h - \left( \eta_h, \mathrm{E}_{h,b} \right)_{L_\epsilon^2 \Lambda^1(\Omega)} ,
    \quad \forall \eta_h \in X_{h,0}^1.
\end{equation*}
With some abuse of notation, we denote now by $\mathbf{M}^1_\epsilon$ the mass matrix corresponding to the space $X^1_h$ (instead of $X^1_{h,0}$). Using the same $0$ and $b$ subindices for the blocks corresponding to $\mathrm{E}_{h,0}$ and $\mathrm{E}_{h,b}$, the mass matrix has the block structure
\[
\mathbf{M}^1_\epsilon = 
\begin{pmatrix}
\mathbf{M}^1_{\epsilon,00} & \mathbf{M}^1_{\epsilon,0b} \\
\mathbf{M}^1_{\epsilon,b0} & \mathbf{M}^1_{\epsilon,bb}
\end{pmatrix},
\]
and with the same subindex notation for the vectors of degrees of freedom, the application of the discrete Hodge star operator for the first scheme becomes
\begin{equation*}
{\bf M}^1_{\epsilon,00} \mathbf{e}_0 = \tilde{\bf K}_{2} \mathbf{d} - {\bf M}^1_{\epsilon,0b} \mathbf{e}_b.
\end{equation*}

{The} application of the boundary condition for the second scheme is done in an analogous way. Once the lifting $\mathrm{E}_{h,b}$ has been computed, the second discrete Hodge star operator is defined as
\begin{equation*}
    \int_\Omega \eta_h \wedge \mathrm{E}_{h,0} = \left( \eta_h, \mathrm{D}_h \right)_{L_{1/\epsilon}^2 \Lambda^2(\Omega)} - \int_\Omega \eta_h \wedge \mathrm{E}_{h,b},
    \quad \forall \eta_h \in \XhT{2}.
\end{equation*}
With a similar abuse of notation as above, we denote now by $\mathbf{K}_1$ the pairing matrix between the spaces $X_h^1$ and $\tilde{X}_h^2$. Obviously, this matrix is not square, but if we split it into blocks corresponding to internal and boundary basis functions as
\[
\mathbf{K}_1 = 
\begin{pmatrix}
(\mathbf{K}_1)_0 & (\mathbf{K}_1)_b
\end{pmatrix}
,
\]
the block $(\mathbf{K}_1)_0$ corresponding to basis functions of $X^1_{h,0}$ is square and invertible. Therefore, the application of the discrete Hodge star operator can be computed by solving 
\begin{equation*}
(\mathbf{K}_1)_0 \mathbf{e}_0 = \tilde{\mathbf{M}}^2_{1/\epsilon} \mathbf{d} - (\mathbf{K}_1)_b \mathbf{e}_b,
\end{equation*}
again for a square matrix.

In both schemes, the boundary condition for $\mathrm{E}_h$ must be also taken into account in the other equations. In particular, the magnetic induction $\mathrm{B}_h$ belongs to $X^2_h$ (instead of $X^2_{h,0}$), and therefore the incidence matrix $\mathbf{D}^1$ must relate the spaces $X^1_h$ and $X^2_h$, adding the rows and columns of boundary functions. Similarly, the matrix $\mathbf{K}_2$ in \eqref{eq:discrete1_hodge_mu} and the matrix $\mathbf{M}^2_{1/\mu}$ in \eqref{eq:discrete2_hodge_mu} must consider, in their columns, the basis functions for the whole space $X^2_h$ instead of $X^2_{h,0}$. Although these matrices are not square anymore, this does not pose any problem, because they are only applied and never inverted.

\section{Fast inversion of pairing matrices with Kronecker product structure}\label{sec:tensorized}
As remarked in the previous section, since pairing matrices only involve the integrals of wedge products, they do not depend on the metric properties of the problem, neither in the form of materials, nor in the geometry of the domain $\Omega$. As a consequence, they can be computed in the parametric domain $(0,1)^3$. We present in this section how the tensor-product structure of the spline spaces can be exploited for the solution of the linear systems involving pairing matrices, and we analyse the computational complexity of the second discretization scheme, given by \eqref{eq:discrete2_ampere}--\eqref{eq:discrete2_hodge_mu}, to provide quantitative results on what kind of gain is achieved with respect to the scheme in which the mass matrix is inverted in \eqref{eq:discrete1_ampere}--\eqref{eq:discrete1_hodge_mu}, closer to canonical approaches available in the literature.

With some abuse of notation, from here onwards we will denote by $S_p(\Xi)$ the univariate space after removing the two boundary functions. The same notation will transfer to the tensor-product spaces.

\subsection{Structure of the pairing matrices}
Let us explore in more detail what the aforementioned tensor-product structure implies. 
To simplify the burden of notation let us initially assume that the knot vectors defining the univariate spline spaces for the space $X_{h,0}^0$ are all the same and denote them as done previously with $\Xi$. Since the wedge product commutes with the pullback, the pairing matrices can be computed in the parametric three-dimensional domain, and the orthogonality of the different Cartesian directions implies that $\mathbf{K}_1$ and $\tilde{\mathbf{K}}_1$ are block diagonal, of the form
\begin{equation*}
\mathbf{K}_1 = 
\begin{bmatrix}
\mathbf{C}_1 & \mathbf{0} & \mathbf{0} \\
\mathbf{0} & \mathbf{C}_2 & \mathbf{0} \\
\mathbf{0} & \mathbf{0} & \mathbf{C}_3 
\end{bmatrix},
\;  \; \;
\tilde{\mathbf{K}}_1 = 
\begin{bmatrix}
\tilde{\mathbf{C}}_1 & \mathbf{0} & \mathbf{0} \\
\mathbf{0} & \tilde{\mathbf{C}}_2 & \mathbf{0} \\
\mathbf{0} & \mathbf{0} & \tilde{\mathbf{C}}_3 
\end{bmatrix}
,
\end{equation*}
with yet to be defined blocks $\mathbf{C}_i$ and $\tilde{\mathbf{C}}_i$. For instance, in the matrix $\mathbf{K}_1$, which gives the pairing between $X_{h,0}^1$ and $\XhT{2}$, the three blocks $\mathbf{C}_i$ are respectively the pairing matrices between the following pairs of discrete spline spaces:
\begin{align*}
  &
  S_{p-1,p,p}(\Xi',\Xi,\Xi) \text{ and } S_{p-1,p-2,p-2}(\Xi',\Xi'',\Xi''),\\
  &
  S_{p,p-1,p}(\Xi,\Xi',\Xi) \text{ and } S_{p-2,p-1,p-2}(\Xi'',\Xi',\Xi''),\\
  &
  S_{p,p,p-1}(\Xi,\Xi,\Xi') \text{ and } S_{p-2,p-2,p-1}(\Xi'',\Xi'',\Xi').
\end{align*}

As a consequence, for the inversion of $\mathbf{K}_1$ (or $\tilde{\mathbf{K}}_1$), the three blocks can be inverted separately and in parallel.
The rest of the section is accordingly devoted to solutions of systems of the type
\begin{equation} \label{eq:diag_systems}
    \mathbf{C}_i \mathbf{x}_i = \mathbf{b}_i, \; \text{ for } i\in\{1,2,3\},
\end{equation}
\noindent in which vectors $\mathbf{x}_i$ and $\mathbf{b}_i$ arise again from the trivial procedure of isolating Cartesian components of the discrete differential forms unknowns.
Focusing on each single block, a further and even more important advantage comes from the fact that each block can be computed by Kronecker tensor product.
To illustrate this let us continue analyzing in detail the matrix $\mathbf{C}_1$. From the definition of the spaces in the parametric domain, this is nothing but the Gram matrix between the basis functions of the tensor-product spaces: 
its entries can then be computed by Kronecker tensor product of the pairing matrices for univariate spaces, since on the parametric domain the integral of the products of univariate basis splines becomes separable. 
Indeed, let $\hat{\mathbf{G}}$ be the pairing matrix between the space of the dual complex $S_{p-2}(\Xi'')$ (rows) and the space of the primal complex $S_{p}(\Xi)$ (columns), and let $\hat{\mathbf{M}}$ be the pairing matrix between the space of the dual complex $S_{p-1}(\Xi')$ (rows) and the space of the primal complex $S_{p-1}(\Xi')$ (columns). Then, we have that 
\begin{equation*}
\mathbf{C}_1 = \hat{\mathbf{G}}\otimes\hat{\mathbf{G}}\otimes\hat{\mathbf{M}}.
\end{equation*}
where the overhead hats stress the fact that entries are computed as integrals on the parametric domain. 
The matrix $\hat{\mathbf{M}}$ is clearly invertible, while $\hat{\mathbf{G}}$ is invertible thanks to the inf-sup condition between the univariate spline spaces of degree $p$ and $p-2$, see \cite{Brivadis_2015aa,Buffa_2020aa}.

We remark that, although pairing the univariate space $S_{p-1}(\Xi)$ with itself, the matrix $\hat{\mathbf{M}}$ is not necessarily symmetric, because we are using two different bases. Concretely, we make the common choice of using the scaled Curry--Schoenberg splines for the univariate spaces of derivatives in the primal complex, i.e. the spaces $S_{p-1}(\Xi')$. For the same univariate space needed in the dual complex we use instead the standard B-spline basis, whereas $S_{p-2}(\Xi'')$ spaces is expanded in the Curry--Schoenberg basis \cite{ratnani}. 
As a consequence, each one dimensional pairing matrix involves a Curry--Schoenberg basis and a standard B-spline basis. Of course, it is possible to decide to use standard B-splines for all involved univariate spaces to recover symmetry, but in this case it is necessary to apply a suitable scaling to the matrices of the exterior derivatives $\mathbf{D}^k$, which would not be incidence matrices anymore, although maintaining the same sparsity pattern.

If we now relax the assumption that the starting polynomial degree $p$ and the knot vector $\Xi$ are the same in every direction, and instead work with different degrees $p_i$ and knot vectors $\Xi_i$, for $i=1,2,3$ no drastic changes occur. The three blocks of the matrix $\mathbf{K}_1$ would take the form
\begin{equation} \label{eq:tensor_K}
\begin{aligned}
\mathbf{C}_1 = \hat{\mathbf{G}}^3 \otimes \hat{\mathbf{G}}^2 \otimes \hat{\mathbf{M}}^1, \\
\mathbf{C}_2 = \hat{\mathbf{G}}^3 \otimes \hat{\mathbf{M}}^2 \otimes \hat{\mathbf{G}}^1, \\
\mathbf{C}_3 = \hat{\mathbf{M}}^3 \otimes \hat{\mathbf{G}}^2 \otimes \hat{\mathbf{G}}^1,
\end{aligned}
\end{equation}
with $\hat{\mathbf{G}}^i$ pairing $S_{p_i-2}(\Xi_i'')$ and $S_{p_i}(\Xi_i)$, and $\hat{\mathbf{M}}^i$ pairing $S_{p_i-1}(\Xi'_i)$ from the dual complex with $S_{p_i-1}(\Xi'_i)$ from the primal complex.

Applying the same reasoning to the blocks of the matrix $\tilde{\mathbf{K}}_1$ leads to very similar expressions. Noting that this is the pairing matrix between the discrete spaces $X_{h,0}^2$ and $\XhT{1}$, and for instance that the first component of these spaces is respectively given by $S_{p,p-1,p-1}(\Xi,\Xi',\Xi')$ and $S_{p-2,p-1,p-1}(\Xi'',\Xi',\Xi')$, the blocks are given by
\begin{equation} \label{eq:tensor_Ktilde}
\begin{aligned}
\tilde{\mathbf{C}}_1 = (\hat{\mathbf{M}}^3)^\top \otimes (\hat{\mathbf{M}}^2)^\top \otimes (\hat{\mathbf{G}}^1)^\top, \\
\tilde{\mathbf{C}}_2 = (\hat{\mathbf{M}}^3)^\top \otimes (\hat{\mathbf{G}}^2)^\top \otimes (\hat{\mathbf{M}}^1)^\top, \\
\tilde{\mathbf{C}}_3 = (\hat{\mathbf{G}}^3)^\top \otimes (\hat{\mathbf{M}}^2)^\top \otimes (\hat{\mathbf{M}}^1)^\top,
\end{aligned}
\end{equation}
where the transposes occur because for $\tilde{\mathbf{K}}_1 = \mathbf{K}_2^\top$, the rows and columns respectively correspond to spaces of the primal and dual complex{, which is the opposite with respect to $\mathbf{K}_1$}.

\subsection{Solution of the linear system}
If a single step time integration scheme is employed, such as the symplectic leapfrog integrator which is common practice for this application, we need to solve a total of six linear systems (three for each Hodge star operator application) where the system matrix is one of the blocks of \eqref{eq:tensor_K} or \eqref{eq:tensor_Ktilde}. Without loss of generality, let us first focus again on matrix $\mathbf{C}_1$, whose inverse is given via a well known property of Kronecker product:
\begin{equation}\label{eq:kroninverse}
  (\hat{\mathbf{G}}^3 \otimes \hat{\mathbf{G}}^2 \otimes \hat{\mathbf{M}}^1)^{-1} = 
(\hat{\mathbf{G}}^3)^{-1} \otimes (\hat{\mathbf{G}}^2)^{-1} \otimes (\hat{\mathbf{M}}^1)^{-1}.
\end{equation}

Following \cite{sangalliIsogeometricPreconditionersBased2016}, for any matrix $\mathbf{X}\in\mathbb{R}^{n_1\times n_2}$ we denote by $\mathrm{vec}(\mathbf{X}) \in \mathbb{R}^{n_1 n_2}$ the vector obtained by unrolling the columns of $\mathbf{X}$ into a single column vector $\mathbf{x}$. Then, if $\mathbf{A}_1$, $\mathbf{A}_2$ are matrices of appropriate dimensions, and $\mathbf{x} = \mathrm{vec}(\mathbf{X})$, the property
\begin{equation*}
(\mathbf{A}_2 \otimes \mathbf{A}_1) \mathbf{x} = \mathrm{vec}(\mathbf{A}_1\mathbf{X}\mathbf{A}_2^\top ),
\end{equation*}
holds. 
This property from tensor algebra can be used to efficiently compute matrix-vector products when the matrix has Kronecker product structure. Indeed, it shows that computing $(\mathbf{A}_2 \otimes \mathbf{A}_1) \mathbf{x}$ is equivalent to computing $n_1$ matrix-vector products with $\mathbf{A}_2$ and $n_2$ matrix-vector products with $\mathbf{A}_1$, and in particular $(\mathbf{A}_2 \otimes \mathbf{A}_1)$ does not ever need to be explicitly computed and stored. 
Furthermore, if $\mathbf{A}_1$ and $\mathbf{A}_2$ are nonsingular, we have an analogous property for the inverse
\begin{equation*}
(\mathbf{A}_2 \otimes \mathbf{A}_1)^{-1} \mathbf{x} = \mathrm{vec}(\mathbf{A}_1^{-1} \mathbf{X} \mathbf{A}_2^{-\top}), 
\end{equation*}
which shows that the problem of solving a linear system for the matrix $(\mathbf{A}_2\otimes\mathbf{A}_1)$ is equivalent to solving $n_1$ linear systems involving $\mathbf{A}_2$ and $n_2$ linear systems involving $\mathbf{A}_1$.

The above property can be used to efficiently solve the linear systems involving the matrices $\mathbf{C}_i$ or $\tilde{\mathbf{C}}_i$. Focusing again on the matrix $\mathbf{C}_1$, and using the structure of its inverse \eqref{eq:kroninverse}, the solution of the linear system is computed as
\begin{equation} \label{eq:system_kronecker}
\mathbf{x}_1 = (\hat{\mathbf{G}}^3 \otimes \hat{\mathbf{G}}^2 \otimes \hat{\mathbf{M}}^1)^{-1} \mathbf{b}_1 = \mathrm{vec} \left((\hat{\mathbf{M}}^1)^{-1} \mathbf{B}_1 (\hat{\mathbf{G}}^3 \otimes \hat{\mathbf{G}}^2)^{-\top}\right).
\end{equation}
Since the rightmost system requires the solution for as many right-hand sides as the number of rows of $\mathbf{B}_1$ (i.e. the dimension of $S_p(\Xi_1)$), it is not possible to further exploit the tensor-product structure again to avoid computing $(\hat{\mathbf{G}}^3 \otimes \hat{\mathbf{G}}^2)^{-\top}$.

{
\begin{remark}
In case of having different meshes in each parametric direction, it is convenient to reorder the directions in such a way that the first direction gets the finest mesh. This would give more similar sizes to the linear systems to be solved in \eqref{eq:system_kronecker}.
\end{remark}
}
\subsection{Computational complexity} \label{sec:complexity}
The study of the computational complexity must be divided into two parts. There is in fact a set of algebraic operations which can be performed only once at the beginning of the time stepping simulation, e.g. assembling the matrices involved in the discrete Hodge star operators and exterior derivatives. 
This system setup is then followed by the cost of applying the exterior derivatives and the discrete Hodge operators during the explicit time stepping. 

For simplicity, we will assume again the same degree and knot vector in every direction, and denote by $n$ the dimension of the univariate space $S_p(\Xi)$. As we are dealing with asymptotic complexity notation, it straightforwardly follows that all dimensions of univariate spline spaces involved in building the discrete spaces of 1--forms and 2--forms are $\mathcal{O}(n)$. We will focus on the simple, but practically very relevant, case of maximal smoothness across all the knots, which implies that the number of elements in each direction is also $\mathcal{O}(n)$. We will denote by $N$ the number of time steps in a simulation, with the reasonable assumption that $N\gg n$ and $N \gg p$, which evidently shifts the focus on the complexity of the time stepping computation rather than the system setup. 
Finally, all complexities reported are based on computing the number of floating point operations (FLOPs) and therein assuming all multiplications and additions involving simple scalars to take constant time.


Due to the local support of splines, all univariate pairing matrices appearing in the Kronecker products have bandwidth $\mathcal{O}(p)$, while the pairing and mass matrices for three-dimensional spaces have dimension $\mathcal{O}(n^3)$ and bandwidth $\mathcal{O}(p^3)$.\footnote{\textcolor{black}{To be precise, it is {not the bandwidth but} the number of nonzeros per row and column.}} Moreover, the incidence matrices for the exterior derivatives have bandwidth $\mathcal{O}(1)$. If the tensor-product structure is not exploited, the complexity of one time step for both discretization schemes \eqref{eq:discrete1_ampere}--\eqref{eq:discrete1_hodge_mu} and \eqref{eq:discrete2_ampere}--\eqref{eq:discrete2_hodge_mu} is bounded by the factorization of their respective system matrices. The complexity of this factorization amounts to the one of matrix-matrix multiplication and hence results in an $\mathcal{O}(n^6 p^3)$ bound, even when considering that the three blocks of the pairing matrices $\mathbf{K}_1$ and $\tilde{\mathbf{K}}_1$ can be solved in parallel for the second scheme. On the other hand, computing the right-hand side in the application of the discrete Hodge star operators involves a matrix-vector multiplication operation of complexity at most $\mathcal{O}(n^3 p^3)$.

To exploit the Kronecker product structure for the second scheme, one must look at the solution of the system in \eqref{eq:system_kronecker}. We first note that the upper and lower triangular matrices of the LU factorization preserve the banded structure, which is $\mathcal{O}(p)$ for the univariate pairing matrices, and $\mathcal{O}(p^2)$ for the matrix $\hat{\mathbf{G}}^3 \otimes \hat{\mathbf{G}}^2$. Proceeding from right to left, one first has to solve $n$ linear systems for a matrix of size $n^2 \times n^2$. Taking into account the bandwidth, this operation has complexity $\mathcal{O}(n^3p^2)$. Subsequently, one has to solve $n^2$ linear systems for a matrix of size $n \times n$ and bandwidth $\mathcal{O}(p)$, an operation of complexity $\mathcal{O}(n^3 p)$. Summing up, the asymptotic cost of solving the linear system at each time step is $\mathcal{O}(n^3p^2)$.

The computational complexity of all other tasks involved in the proposed algorithm which exploits the Kronecker product structure are summarized in Table~\ref{tab:complexity}. It is then clear that applying the mass matrix is the most consuming task in terms of FLOPs, {even more than solving the linear systems,} while still being linear in the number of unknowns, which in three dimensional problems grow as $n^3$. Therefore, we obtain an estimate of the computational complexity for the second scheme of $\mathcal{O}(N n^3 p^3)$. This is a substantial gain with respect to the scheme in which we have to solve the linear system for the mass matrices, and which has computational complexity of order $\mathcal{O}(N n^6 p^3)$.

\begin{table}[t]
  \renewcommand{\arraystretch}{1.3}
    \centering
\begin{tabular}{|c|c|c|c|c|c|c}\hline
    \multicolumn{7}{|c|}{\textbf{Complexity of tensorized solution of \eqref{eq:discrete2_ampere}--\eqref{eq:discrete2_hodge_mu} explicit in} $n$, $p$, $N$}\\\hline
    Procedure & \multicolumn{4}{c|}{Notes on operation} & \multicolumn{2}{c|}{Complexity} \\\hline
    Assemble $\mathbf{M}_{\mu^{-1}}^2$, $\tilde{\mathbf{M}}_{\epsilon^{-1}}^2$ & \multicolumn{4}{c|}{Numerical integration, banded matrix} & \multicolumn{2}{c|}{$\mathcal{O}(p^9 n^{3})$} \\
    Assemble $\mathbf{D}^1$, $\tilde{\mathbf{D}}^1$ & \multicolumn{4}{c|}{Scaled incidence matrix} & \multicolumn{2}{c|}{$\mathcal{O}(n^3)$} \\
    Assemble $\hat{\mathbf{G}}^k, \hat{\mathbf{M}}^k$, $k=1,2,3$ & \multicolumn{4}{c|}{Numerical integration, univariate spaces} & \multicolumn{2}{c|}{$\mathcal{O}(p^3 n)$} \\
    Factorization of $\hat{\mathbf{G}}^k, \hat{\mathbf{M}}^k$, $k=1,2,3$ & \multicolumn{4}{c|}{Univariate spaces, banded matrices } & \multicolumn{2}{c|}{$\mathcal{O}(p^2 n)$} \\
\hline
    -- & \multicolumn{4}{c|}{Assembly of the whole system} & \multicolumn{2}{c|}{$\mathcal{O}(p^9 n^{3})$} \\ \hline
    Apply $\mathbf{D}^1$, $\tilde{\mathbf{D}}^1$ & \multicolumn{4}{c|}{Banded matrix-vector product} & \multicolumn{2}{c|}{$\mathcal{O}( N n^3)$} \\ 
    Apply $\mathbf{M}_{\mu^{-1}}^2$, $\tilde{\mathbf{M}}_{\epsilon^{-1}}^2$ & \multicolumn{4}{c|}{Banded matrix-vector product} & \multicolumn{2}{c|}{$\mathcal{O}( N n^3 p^3)$} \\ 
    Solve the systems for $\mathbf{K}_{1}$, $\tilde{\mathbf{K}}_{1}$  & \multicolumn{4}{c|}{Banded matrix-matrix product} & \multicolumn{2}{c|}{$\mathcal{O}(N n^3 p^2)$} \\ \hline
    -- & \multicolumn{4}{c|}{Time stepping} & \multicolumn{2}{c|}{$\mathcal{O}( N n^3 p^3)$} \\ \hline
    \end{tabular}
    \label{tab:complexity}
    \caption{All major linear algebra computation involved in the time dependent simulation with their respective computation complexity, where $n$ is the dimension of the univariate B-spline space, $p$ is the polynomial degree, and $N$ is the number of time steps.}
\end{table}

\begin{remark}
The Kronecker structure of the pairing matrices can be also exploited in the first scheme, for an efficient computation of matrix-vector products in the application of the Hodge operators. However, the dominant part will remain the solution of the linear systems for the mass matrices.
\end{remark}

\begin{remark}
The incidence matrices $\mathbf{D}^k$ are also independent of the metric properties of the domain, and they can be computed by Kronecker products of univariate incidence matrices and identity matrices, see for instance \cite{Holderied_2021}. Since their application is far from the bottleneck of computation, the details regarding their optimization are neglected.
\end{remark}

\section{Numerical results}\label{sec:num}
In the present section we discuss the behavior of the proposed method, in terms of accuracy, conservation of energy and computational cost, on several numerical tests. In the following we will label all tests by their starting polynomial degree $p$ for 0--forms in the primal discrete complex, since all other degrees follow by sequence properties. Furthermore, the regularity of splines will always be the maximum available one, which is $p-1$ for univariate splines of degree $p$. For time integration we will apply a low order leapfrog method. All tests are carried out by implementing the underlying space and time integration numerical schemes in MATLAB within the open-source library GeoPDEs~\cite{Vazquez_2016aa}, an IGA based software library oriented towards basic research. 

\subsection{Accuracy of the spatial discretization and conservation of energy}
The first results concern the convergence of the method with respect to the mesh size, and the conservation of energy.

\subsubsection{Cavity problem in the unit cube}
The first geometry used as test example is the unit cube (i.e. $\Omega=\hat{\Omega}=[0,1]^3$) with homogeneous boundary conditions on the electric field 1--form, which implies studying a Maxwell cavity problem with a known time harmonic solution {(the first eigenfunction of the double curl operator for the cube) with time period $\sqrt{2} \pi$.} The physical parameters are all normalized such that the dielectric permittivity and magnetic permeability of the vacuum are $\epsilon_0=\mu_0=1$. The initial conditions for $\mathrm{B}$ and $\mathrm{D}$ are set by projecting their value for $t=0$ into the discrete 2--form spaces.

\begin{figure}[!t]
  \centering
  \begin{tikzpicture}
    \begin{axis}[width=0.9\textwidth,  height=0.35\textwidth,
       ylabel={$\mathcal{E}_h(t)$},
      xlabel={Time [natural units]},
          ymin=0.28,ymax=0.4,
        xmin=0.0,xmax=20, 
        ylabel near ticks,
        tick label style={
            /pgf/number format/.cd,
                fixed relative,
                precision=3,
                zerofill,
            /tikz/.cd,
        },
        legend columns=3, 
        legend style={
            /tikz/column 2/.style={
                column sep=5pt,
            },
        },
        label style={font=\normalsize},every x tick scale label/.style={at={(1,0)},
        anchor=north,yshift=-5pt,inner sep=0pt},legend pos= south east]
    \addplot [red,dashed,thick] table [x index=0, y index=2, col sep=comma] {figures/energy_conservation/maxwell_3d_cube_performance_decompose_K_tp_jcp_p3r2_nsubs2_subs_energy.csv};
    \addplot [black,dotted,thick] table [x index=0, y index=2, col sep=comma] {figures/energy_conservation/maxwell_3d_cube_performance_decompose_K_tp_jcp_p3r2_nsubs4_subs_energy.csv};
    \addplot [blue,solid,thick] table [x index=0, y index=2, col sep=comma] {figures/energy_conservation/maxwell_3d_cube_performance_decompose_K_tp_jcp_p3r2_nsubs8_subs_energy.csv};
    \legend{  
    $p=3$, $h=0.5$\\
    $p=3$, $h=0.25$\\
    $p=3$, $h=0.125$\\
    }
\end{axis}
  \end{tikzpicture}\\
  \begin{tikzpicture}
    \begin{axis}[width=0.9\textwidth,  height=0.35\textwidth,
       ylabel={$\mathcal{E}_h(t)$},
      xlabel={Time [natural units]},
      ymin=0.28,ymax=0.4,
        xmin=0.0,xmax=20, 
        ylabel near ticks,
        tick label style={
            /pgf/number format/.cd,
                fixed relative,
                precision=3,
                zerofill,
            /tikz/.cd,
        },
        legend columns=3, 
        legend style={
            /tikz/column 2/.style={
                column sep=5pt,
            },
        },
        label style={font=\normalsize},every x tick scale label/.style={at={(1,0)},
        anchor=north,yshift=-5pt,inner sep=0pt},legend pos= south east]
    \addplot [red,dashed,thick] table [x index=0, y index=2, col sep=comma] {figures/energy_conservation/maxwell_3d_cube_performance_decompose_K_tp_jcp_p4r3_nsubs2_subs_energy.csv};
    \addplot [black,dotted,thick] table [x index=0, y index=2, col sep=comma] {figures/energy_conservation/maxwell_3d_cube_performance_decompose_K_tp_jcp_p4r3_nsubs4_subs_energy.csv};
    \addplot [blue,solid,thick] table [x index=0, y index=2, col sep=comma] {figures/energy_conservation/maxwell_3d_cube_performance_decompose_K_tp_jcp_p4r3_nsubs8_subs_energy.csv};
    \legend{  
    $p=4$, $h=0.5$\\
    $p=4$, $h=0.25$\\
    $p=4$, $h=0.125$\\
    }
\end{axis}
  \end{tikzpicture}
  \caption{ 
    The discrete energy over time for the approach \eqref{eq:discrete2_ampere}--\eqref{eq:discrete2_hodge_mu} for degrees $p=3,4$ and different mesh sizes.
    } \label{fig:energy_vs_time}
\end{figure}
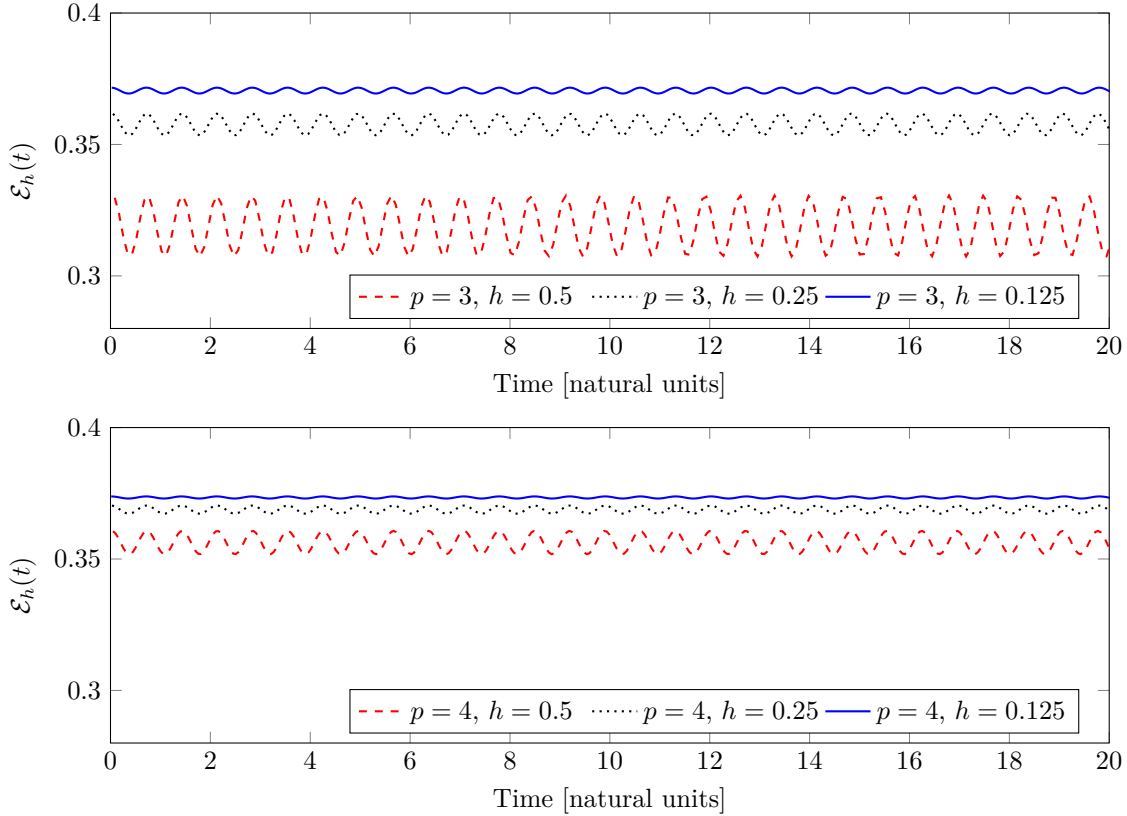

The cavity problem is used to test the proven energy conservation properties of the proposed method, both under $p$ and $h$ refinement (with $h$ being the spatial mesh size), and in particular we solve the problem for degrees $p=3,4$ on successively refined meshes, with the time step given by the Courant-Friedrichs-Lewy (CFL) condition. The results shown in Fig.~\ref{fig:energy_vs_time}, which correspond to the second discretization scheme of Section~\ref{sec:solve_K}, confirm the theoretical result about the conservation of energy, {with oscillations appearing due to the typical behavior of the leapfrog scheme. The amplitude of the oscillations is reduced for higher degree and finer meshes, and also the error with respect to the exact energy reduces both when refining the mesh and when increasing the degree.}

\begin{figure}[t]
\begin{tikzpicture}[]
  \begin{loglogaxis}[width=0.46\textwidth, height=0.46\textwidth,
    ylabel={$\|{\mathrm{E}}_h - {\mathrm{E}}\|_{L^2} / \|\mathrm{E}\|_{L^2}$},xlabel={Mesh size}, 
    ymin=1E-6,ymax=2E2,
      xmin=5E-2,xmax=1, 
      ylabel near ticks,tick label style={font=\normalsize},    
                    legend columns=2, 
              legend style={font=\scriptsize,
                  /tikz/column 2/.style={
                      column sep=5pt,
                  },
              },
      label style={font=\normalsize},every x tick scale label/.style={at={(1,0)},
      anchor=north,yshift=-5pt,inner sep=0pt},legend pos= north west]
  \addplot [blue,solid,thick,mark=triangle*, mark size=3pt] table [x index=0, y index=11, col sep=comma] {figures/relative_errors/cube/p2/maxwell_3d_cube_relative_invertK_tp_jcp_p2r1_metrics.csv};
  \addplot [red,solid,thick,mark=triangle, mark size=3pt]  table [x index=0, y index=11, col sep=comma] {figures/relative_errors/cube/p2/maxwell_3d_cube_relative_invertM_jcp_p2r1_metrics.csv};
  \addplot [blue,solid,thick,mark=square*, mark size=2.5pt] table [x index=0, y index=11, col sep=comma] {figures/relative_errors/cube/p3/maxwell_3d_cube_relative_invertK_tp_jcp_p3r2_metrics.csv};
  \addplot [red,solid,thick,mark=square, mark size=2.5pt]  table [x index=0, y index=11, col sep=comma] {figures/relative_errors/cube/p3/maxwell_3d_cube_relative_invertM_jcp_p3r2_metrics.csv};
  \addplot [blue,solid,thick,mark=*, mark size=2.5pt] table [x index=0, y index=11, col sep=comma] {figures/relative_errors/cube/p4/maxwell_3d_cube_relative_invertK_tp_jcp_p4r3_metrics.csv};
  \addplot [red,solid,thick,mark=o, mark size=2.5pt]  table [x index=0, y index=11, col sep=comma] {figures/relative_errors/cube/p4/maxwell_3d_cube_relative_invertM_jcp_p4r3_metrics.csv};
  
  \convergenceslopeinvh{7e-2}{12e-2}{0.1}{2}{}{};
  \convergenceslopeinvh{7e-2}{12e-2}{2E-4}{3}{}{};
  \convergenceslopeh{2e-1}{3e-1}{1E-4}{4}{}{};
  \legend{  
    Solve $\mathbf{K}$, $p=2$\\
    Solve $\mathbf{M}$, $p=2$\\
    Solve $\mathbf{K}$, $p=3$\\
    Solve $\mathbf{M}$, $p=3$\\
    Solve $\mathbf{K}$, $p=4$\\
    Solve $\mathbf{M}$, $p=4$\\}
  \end{loglogaxis}
\end{tikzpicture}
\begin{tikzpicture}[]
  \begin{loglogaxis}[width=0.46\textwidth, height=0.46\textwidth,
    ylabel={$\|{\mathrm{H}}_h - {\mathrm{H}}\|_{L^2} / \|\mathrm{H}\|_{L^2}$},xlabel={Mesh size}, 
    ymin=1E-6,ymax=2E2,
      xmin=5E-2,xmax=1, 
      ylabel near ticks,tick label style={font=\normalsize},  
                    legend columns=2, 
              legend style={font=\scriptsize,
                  /tikz/column 2/.style={
                      column sep=5pt,
                  },
              },
      label style={font=\normalsize},every x tick scale label/.style={at={(1,0)},
      anchor=north,yshift=-5pt,inner sep=0pt},legend pos= north east]
  \addplot [blue,solid,thick,mark=triangle*, mark size=3pt] table [x index=0, y index=13, col sep=comma] {figures/relative_errors/cube/p2/maxwell_3d_cube_relative_invertK_tp_jcp_p2r1_metrics.csv};
  \addplot [red,solid,thick,mark=triangle, mark size=3pt]  table [x index=0, y index=13, col sep=comma] {figures/relative_errors/cube/p2/maxwell_3d_cube_relative_invertM_jcp_p2r1_metrics.csv};
  \addplot [blue,solid,thick,mark=square*, mark size=2.5pt] table [x index=0, y index=13, col sep=comma] {figures/relative_errors/cube/p3/maxwell_3d_cube_relative_invertK_tp_jcp_p3r2_metrics.csv};
  \addplot [red,solid,thick,mark=square, mark size=2.5pt]  table [x index=0, y index=13, col sep=comma] {figures/relative_errors/cube/p3/maxwell_3d_cube_relative_invertM_jcp_p3r2_metrics.csv};
  \addplot [blue,solid,thick,mark=*, mark size=2.5pt] table [x index=0, y index=13, col sep=comma] {figures/relative_errors/cube/p4/maxwell_3d_cube_relative_invertK_tp_jcp_p4r3_metrics.csv};
  \addplot [red,solid,thick,mark=o, mark size=2.5pt]  table [x index=0, y index=13, col sep=comma] {figures/relative_errors/cube/p4/maxwell_3d_cube_relative_invertM_jcp_p4r3_metrics.csv};
  \convergenceslopeinvh{7e-2}{12e-2}{0.2}{1}{}{};
  \convergenceslopeinvh{7e-2}{12e-2}{9E-3}{2}{}{};
  \convergenceslopeinvh{7e-2}{12e-2}{6E-4}{3}{}{};
  \legend{  
    Solve $\mathbf{K}$, $p=2$\\
    Solve $\mathbf{M}$, $p=2$\\
    Solve $\mathbf{K}$, $p=3$\\
    Solve $\mathbf{M}$, $p=3$\\
    Solve $\mathbf{K}$, $p=4$\\
    Solve $\mathbf{M}$, $p=4$\\}
  \end{loglogaxis}
\end{tikzpicture}
\caption{Relative error convergence rates in the $L^2(0,T;L^2(\Omega))$ norm for the unit cube, for the electric field (left) and the magnetic field (right) with polynomial degrees $p=2,3,4$. We use label $\mathbf{M}$ to label the first discretization scheme, and $\mathbf{K}$ for the second discretization scheme. } \label{fig:convergence_cube}
\end{figure}
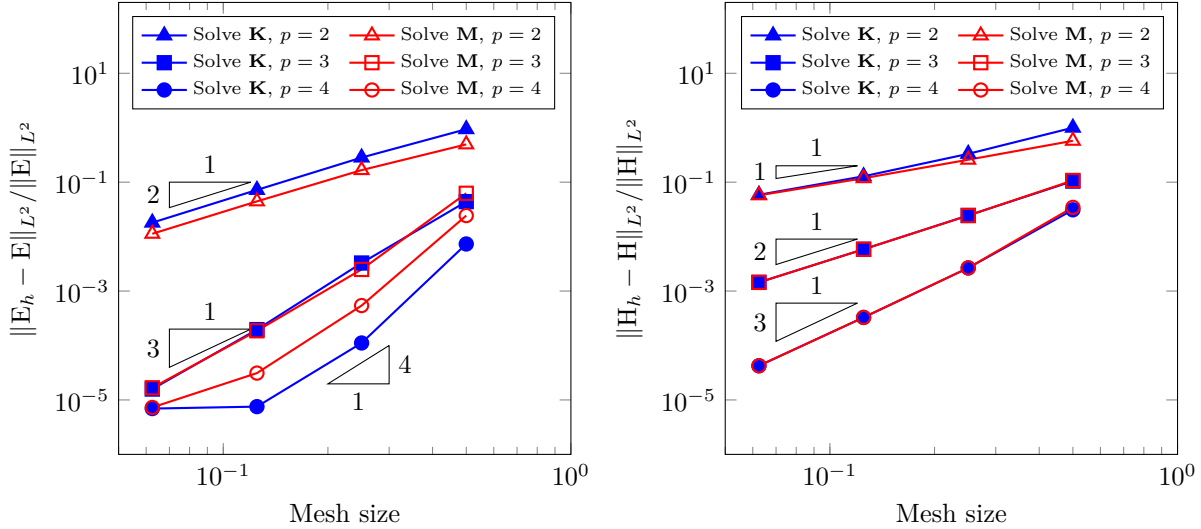

To test the order of convergence of the two methods under mesh refinement we solve the problem in the time interval {$[0,T]$ with $T=2.0$.} We consider now degrees $p=2,3,4$, and different meshes refined dyadically, with the number of elements ranging from $2$ to $16$ in each Cartesian direction. In order to observe the spatial discretization error, a very fine step of value $\Delta t_{\text{ref}}=8.0548\times 10^{-4}$, measured in natural units, is used for all degrees and meshes. This value is prescribed by the computed CFL condition for the most refined mesh among all numerical tests {(including the coaxial cable below)}, which is then divided further by a factor ten. Since we are approximating the electric field $\mathrm{E}$ and the magnetic field $\mathrm{H}$ in two different sequences, we compute the {relative} error for both of them, using the $L^2(0,T;L^2(\Omega))$ norm, that is, $L^2$ norm both in space and time. {Although not reported, similar results are obtained with the $L^\infty$ norm in time.} The convergence plots are shown in Fig.~\ref{fig:convergence_cube}, where we observe that with both schemes the error converges to zero as $h^p$ for the electric field, which is approximated with spaces of mixed degree $p$ and $p-1$, while for the magnetic field, which is approximated with mixed degrees $p-1$ and $p-2$, the error converges to zero as $h^{p-1}$. This error is consistent with the ones that would be obtained by approximation with finite element methods of the same degrees. Moreover, the magnitude of the error is very similar for both schemes. We note that a plateau is observed for the electric field with $p=4$ at the finest level, caused by the time discretization error. This effect could be removed by using a time discretization scheme of higher order{, or a finer time step.}

\subsubsection{Coaxial cable}
\begin{figure}[!t]
  \centering
  \includegraphics[width=0.49\textwidth]{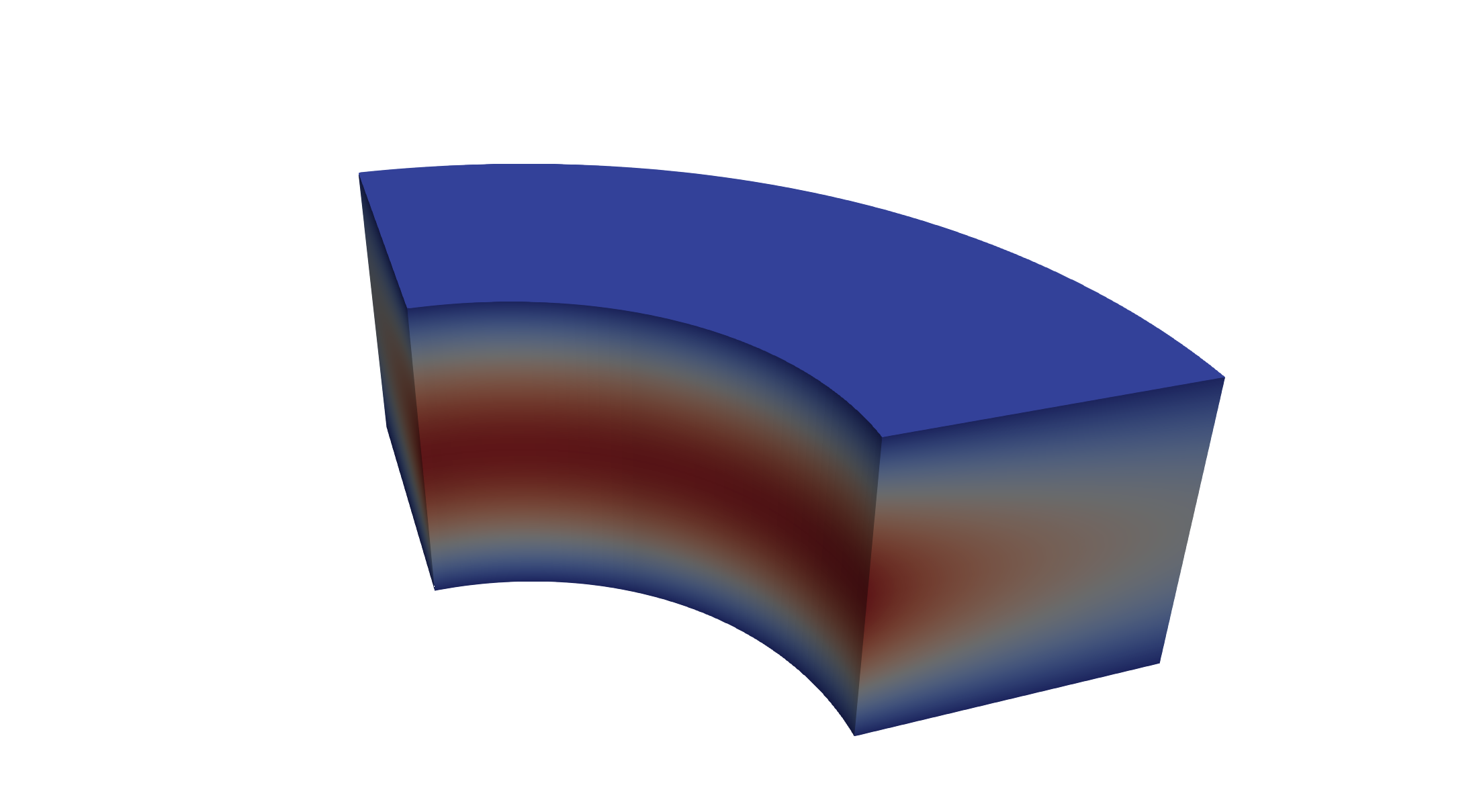}
  \includegraphics[width=0.49\textwidth]{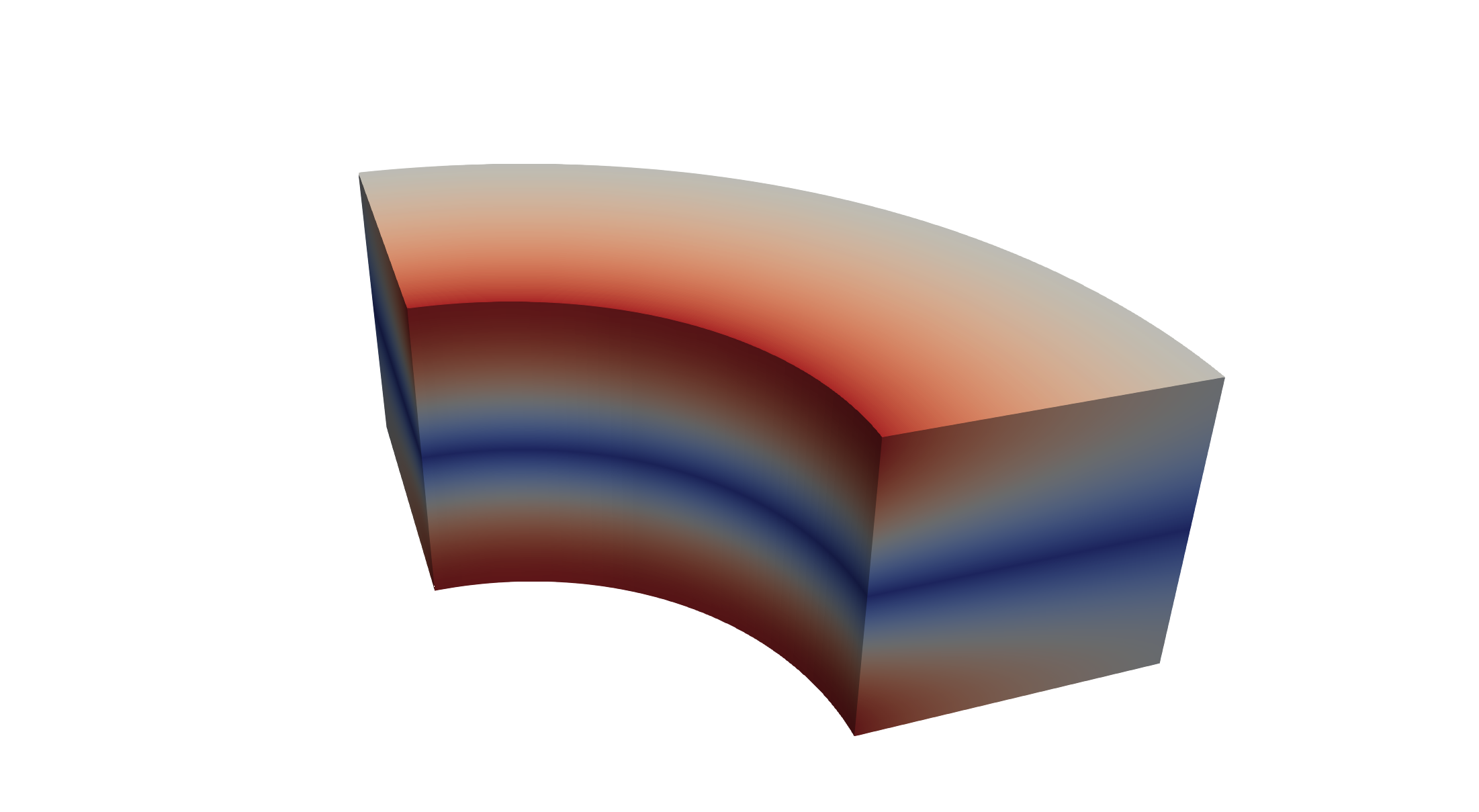}
  \caption{A quarter of a coaxial cable and the magnitude of the electric field (left) and the magnetic field (right) for its TEM mode.}
  \label{fig:wg_bend}
\end{figure}

For the second numerical test we consider the geometry of a coaxial cable, for which we simulate the steady state propagation of the transverse-electro-magnetic (TEM) mode, the solution in polar coordinates can be found in \cite[Chapter~3]{collinFoundationsMicrowaveEngineering2001}. We just model one quarter of the full coaxial cable $\Omega = \{(x,y,z) : 1 < x^2 + y^2 < 2, 0<z<1 \}$, described with a NURBS geometry of degree 2, and exploit symmetries in the solution, shown in Fig.~\ref{fig:wg_bend} for $t=0.5$, computed with starting degree $p=3$ and maximum regularity in the primal sequence and a mesh of eight elements in each parametric direction, and the same reference time step as mentioned above. As above, we normalize the values of $\epsilon_0$ and $\mu_0$, and set the initial conditions for $\mathrm{B}$ and $\mathrm{D}$ by projecting the exact initial condition into the discrete spaces of 2--forms. Moreover, at every time step we impose inhomogeneous Dirichlet boundary conditions (for $y=0$ and $x=0$), as discussed in Section~\ref{sec:hodge}. As the boundary condition is given by a separable function of space and time, its contribution to the right-hand side is performed only once at the start of numerical time integration, and then the representation in terms of degrees of freedom is multiplied by a known function of time at each time step.

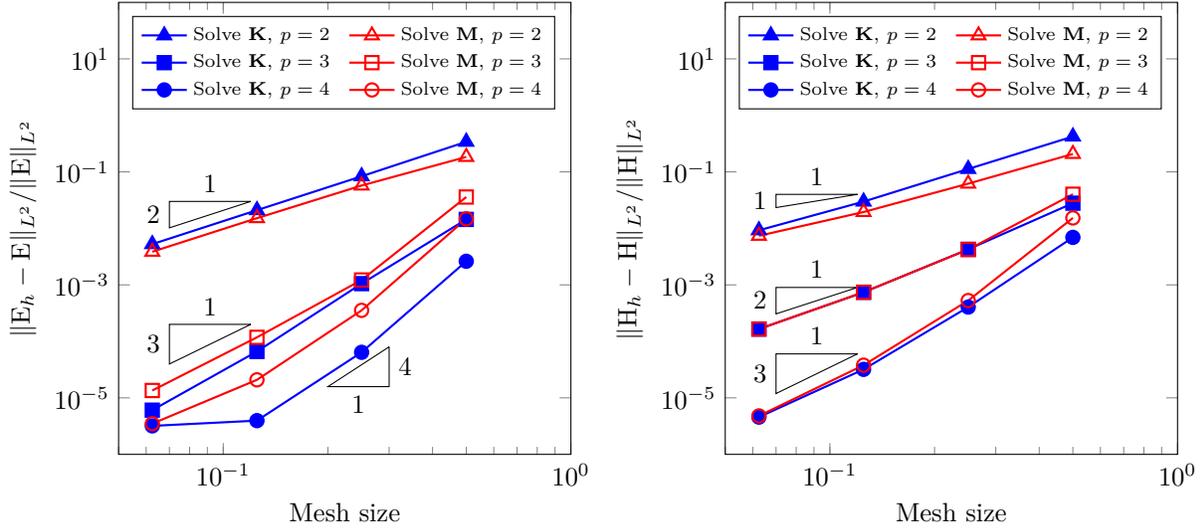
\begin{figure}[t]
\begin{tikzpicture}[]
  \begin{loglogaxis}[width=0.46\textwidth, height=0.46\textwidth,
    ylabel={$\|{\mathrm{E}}_h - {\mathrm{E}}\|_{L^2} / \|\mathrm{E}\|_{L^2}$},xlabel={Mesh size}, 
    ymin=1E-6,ymax=1E2,
      xmin=5E-2,xmax=1, 
      ylabel near ticks,tick label style={font=\normalsize},    
                    legend columns=2, 
              legend style={font=\scriptsize,
                  /tikz/column 2/.style={
                      column sep=5pt,
                  },
              },
      label style={font=\normalsize},every x tick scale label/.style={at={(1,0)},
      anchor=north,yshift=-5pt,inner sep=0pt},legend pos= north west]
  \addplot [blue,solid,thick,mark=triangle*, mark size=3pt] table [x index=0, y index=11, col sep=comma] {figures/relative_errors/coaxial_cable/p2/maxwell_3d_coaxial_relative_invertK_tp_jcp_p2r1_metrics.csv};
  \addplot [red,solid,thick,mark=triangle, mark size=3pt]  table [x index=0, y index=11, col sep=comma] {figures/relative_errors/coaxial_cable/p2/maxwell_3d_coaxial_relative_invertM_jcp_p2r1_metrics.csv};
  \addplot [blue,solid,thick,mark=square*, mark size=2.5pt] table [x index=0, y index=11, col sep=comma] {figures/relative_errors/coaxial_cable/p3/maxwell_3d_coaxial_relative_invertK_tp_jcp_p3r2_metrics.csv};
  \addplot [red,solid,thick,mark=square, mark size=2.5pt]  table [x index=0, y index=11, col sep=comma] {figures/relative_errors/coaxial_cable/p3/maxwell_3d_coaxial_relative_invertM_jcp_p3r2_metrics.csv};
  \addplot [blue,solid,thick,mark=*, mark size=2.5pt] table [x index=0, y index=11, col sep=comma] {figures/relative_errors/coaxial_cable/p4/maxwell_3d_coaxial_relative_invertK_tp_jcp_p4r3_metrics.csv};
  \addplot [red,solid,thick,mark=o, mark size=2.5pt]  table [x index=0, y index=11, col sep=comma] {figures/relative_errors/coaxial_cable/p4/maxwell_3d_coaxial_relative_invertM_jcp_p4r3_metrics.csv};
  \convergenceslopeinvh{7e-2}{12e-2}{0.03}{2}{}{};
  \convergenceslopeinvh{7e-2}{12e-2}{2E-4}{3}{}{};
  \convergenceslopeh{2e-1}{3e-1}{8E-5}{4}{}{};
  \legend{  
    Solve $\mathbf{K}$, $p=2$\\
    Solve $\mathbf{M}$, $p=2$\\
    Solve $\mathbf{K}$, $p=3$\\
    Solve $\mathbf{M}$, $p=3$\\
    Solve $\mathbf{K}$, $p=4$\\
    Solve $\mathbf{M}$, $p=4$\\}
  \end{loglogaxis}
\end{tikzpicture}
\begin{tikzpicture}[]
  \begin{loglogaxis}[width=0.46\textwidth, height=0.46\textwidth,
    ylabel={$\|{\mathrm{H}}_h - {\mathrm{H}}\|_{L^2} / \|\mathrm{H}\|_{L^2}$},xlabel={Mesh size}, 
    ymin=1E-6,ymax=1E2,
      xmin=5E-2,xmax=1, 
      ylabel near ticks,tick label style={font=\normalsize},    
                    legend columns=2, 
              legend style={font=\scriptsize,
                  /tikz/column 2/.style={
                      column sep=5pt,
                  },
              },
      label style={font=\normalsize},every x tick scale label/.style={at={(1,0)},
      anchor=north,yshift=-5pt,inner sep=0pt},legend pos= north west]
  \addplot [blue,solid,thick,mark=triangle*, mark size=3pt] table [x index=0, y index=13, col sep=comma] {figures/relative_errors/coaxial_cable/p2/maxwell_3d_coaxial_relative_invertK_tp_jcp_p2r1_metrics.csv};
  \addplot [red,solid,thick,mark=triangle, mark size=3pt]  table [x index=0, y index=13, col sep=comma] {figures/relative_errors/coaxial_cable/p2/maxwell_3d_coaxial_relative_invertM_jcp_p2r1_metrics.csv};
  \addplot [blue,solid,thick,mark=square*, mark size=2.5pt] table [x index=0, y index=13, col sep=comma] {figures/relative_errors/coaxial_cable/p3/maxwell_3d_coaxial_relative_invertK_tp_jcp_p3r2_metrics.csv};
  \addplot [red,solid,thick,mark=square, mark size=2.5pt]  table [x index=0, y index=13, col sep=comma] {figures/relative_errors/coaxial_cable/p3/maxwell_3d_coaxial_relative_invertM_jcp_p3r2_metrics.csv};
  \addplot [blue,solid,thick,mark=*, mark size=2.5pt] table [x index=0, y index=13, col sep=comma] {figures/relative_errors/coaxial_cable/p4/maxwell_3d_coaxial_relative_invertK_tp_jcp_p4r3_metrics.csv};
  \addplot [red,solid,thick,mark=o, mark size=2.5pt]  table [x index=0, y index=13, col sep=comma] {figures/relative_errors/coaxial_cable/p4/maxwell_3d_coaxial_relative_invertM_jcp_p4r3_metrics.csv};
  \convergenceslopeinvh{7e-2}{12e-2}{0.04}{1}{}{};
  \convergenceslopeinvh{7e-2}{12e-2}{9E-4}{2}{}{};
  \convergenceslopeinvh{7e-2}{12e-2}{6E-5}{3}{}{};
  \legend{  
    Solve $\mathbf{K}$, $p=2$\\
    Solve $\mathbf{M}$, $p=2$\\
    Solve $\mathbf{K}$, $p=3$\\
    Solve $\mathbf{M}$, $p=3$\\
    Solve $\mathbf{K}$, $p=4$\\
    Solve $\mathbf{M}$, $p=4$\\}
  \end{loglogaxis}
\end{tikzpicture}
\caption{Relative error convergence rates in the $L^2(0,T;L^2(\Omega))$ norm for the coaxial cable, for the electric field (left) and the magnetic field (right) with polynomial degrees $p=2,3,4$. We use label $\mathbf{M}$ to label the first discretization scheme, and $\mathbf{K}$ for the second discretization scheme. } \label{fig:convergence_cable}
\end{figure}

In this numerical test we assess the metric dependence of the Hodge star operators. We run a convergence test for both methods with the same degrees, number of elements, and time step as for the unit cube, and we compute the {relative} errors for the electric and magnetic fields as before, the results are shown in Fig.~\ref{fig:convergence_cable}. We observe that the error converges to zero with the same rates as for the unit cube, in which the metric was not involved, i.e., as $h^p$ for the electric field and as $h^{p-1}$ for the magnetic field. Moreover, the magnitude of the error for the two methods is very similar, without a clear advantage for one or the other approach.

\subsection{Study of the computational cost}
In terms of computational efficiency, we aim at validating the computational complexity estimates of Section~\ref{sec:complexity}, for which we use the same numerical examples as for the study of the accuracy. Since those estimates depend on the number of time steps, we also analyse the behavior of the CFL condition with respect to the mesh size and the degree of the splines.

\begin{figure}[!t]
  \centering
\begin{tikzpicture}[]
    \begin{loglogaxis}[width=0.46\textwidth, height=0.46\textwidth,
      ylabel={$\Delta t$},
      xlabel={Polynomial degree},
      title={$\Delta t_{max}$ versus $p$},
        ymin=0.005,ymax=1,xmin=1,xmax=9, ylabel near ticks,tick label style={font=\normalsize},  legend style={font=\normalsize,nodes={scale=0.8, transform shape}}, label style={font=\normalsize},every x tick scale label/.style={at={(1,0)},anchor=north,yshift=-5pt,inner sep=0pt},legend pos= south west]
    \addplot [blue,thick,mark=*, mark size=2.5pt] table [x index=0, y index=1, col sep=comma] {figures/maxwell_timestep_bend_pref.csv};
    \addplot [red,thick,mark=triangle*, mark size=3pt]      table [x index=0, y index=2, col sep=comma] {figures/maxwell_timestep_bend_pref.csv};
    \convergenceslope{3}{5}{0.1}{2}{}{};
    \convergenceslopeinv{4}{6}{0.18}{1.5}{}{};
    \legend{ Solve $\mathbf{K}$  \\ Solve $\mathbf{M}$ \\}
    \end{loglogaxis}
    \end{tikzpicture}
    \begin{tikzpicture}[]
      \begin{loglogaxis}[width=0.45\textwidth,  height=0.45\textwidth,
        title={$\Delta t_{max}$ versus $h$},
        ylabel={$\Delta t$},
        xlabel={Mesh size},  
          ymin=5E-3,ymax=1,
          xmin=5E-2,xmax=7e-1, 
          ylabel near ticks,tick label style={font=\normalsize},  
          legend style={font=\normalsize,nodes={scale=0.8, transform shape}}, 
          legend columns=2, 
          label style={font=\normalsize},every x tick scale label/.style={at={(1,0)},
          anchor=north,yshift=-5pt,inner sep=0pt},legend pos= north west]
      \addplot [blue,solid,thick,mark=triangle*, mark size=3pt] table [x index=0, y index=3, col sep=comma] {figures/timestep_study/maxwell_3d_onestep_invertK_jcp_p2r1_metrics.csv};
      \addplot [red,dashed,thick,mark=triangle, mark size=3pt, mark options={solid,fill=red}]  table [x index=0, y index=3, col sep=comma] {figures/timestep_study/maxwell_3d_onestep_invertM_jcp_p2r1_metrics.csv};
      \addplot [blue,solid,thick,mark=square*, mark size=2.5pt] table [x index=0, y index=3, col sep=comma] {figures/timestep_study/maxwell_3d_onestep_invertK_jcp_p3r2_metrics.csv};
      \addplot [red,dashed,thick,mark=square, mark size=2.5pt, mark options={solid,fill=red}]  table [x index=0, y index=3, col sep=comma] {figures/timestep_study/maxwell_3d_onestep_invertM_jcp_p3r2_metrics.csv};
      \addplot [blue,solid,thick,mark=*, mark size=2.5pt] table [x index=0, y index=3, col sep=comma] {figures/timestep_study/maxwell_3d_onestep_invertK_jcp_p4r3_metrics.csv};
      \addplot [red,dashed,thick,mark=o, mark size=2.5pt, mark options={solid,fill=red}]  table [x index=0, y index=3, col sep=comma] {figures/timestep_study/maxwell_3d_onestep_invertM_jcp_p4r3_metrics.csv};
      \convergenceslopeh{0.1}{0.2}{0.02}{1}{}{};
      \legend{  
        Solve $\mathbf{K}$, $p=2$ \\
        Solve $\mathbf{M}$, $p=2$ \\
        Solve $\mathbf{K}$, $p=3$ \\
        Solve $\mathbf{M}$, $p=3$ \\
        Solve $\mathbf{K}$, $p=4$ \\
        Solve $\mathbf{M}$, $p=4$ \\}
      \end{loglogaxis}
    \end{tikzpicture}
\caption{The maximum time step given by the CFL condition with respect to the polynomial degree (left), and with respect to the mesh size (right).} 
\label{fig:cfl_vs_p_and_h}
\end{figure}
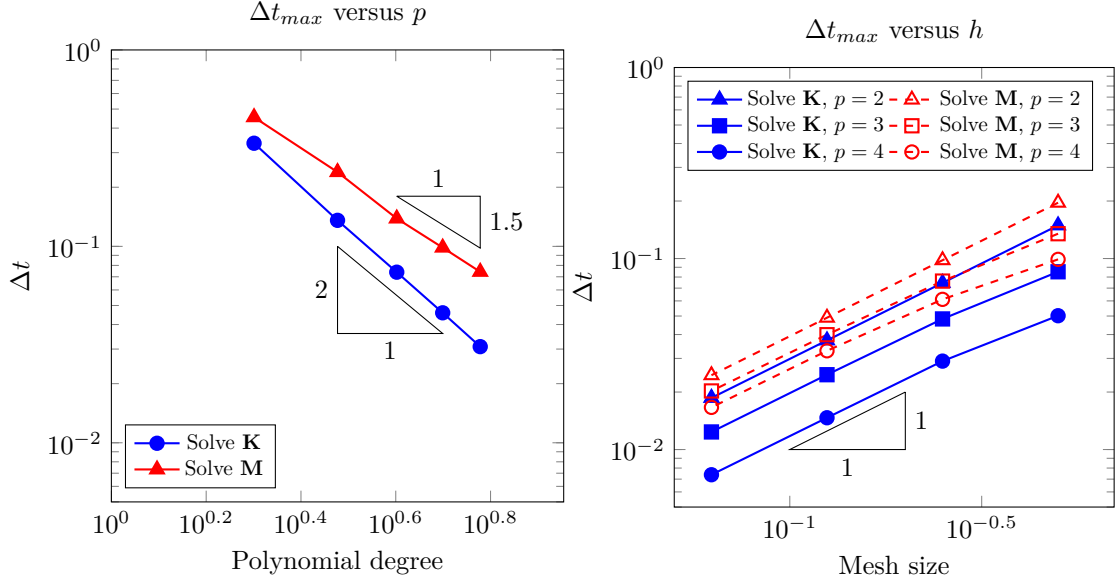

To understand the behavior of the CFL condition, which gives the maximum allowed time step, we compute it for both methods for degrees from 2 to 6, and also for different mesh refinements of the unit cube. The plot on the left of Fig.~\ref{fig:cfl_vs_p_and_h}, in which we vary the degree and fix the number of elements to the one of the finest level in the $h$--refinement, shows that for the second discretization scheme, in which we have to solve the linear systems for the pairing matrix, the CFL condition scales as $\mathcal{O}(1/p^2)$, which is the standard behavior for FEM and DG-FEM. Instead, for the scheme based on solving linear systems associated to the mass matrix, the CFL condition scales as $\mathcal{O}(1/p^{3/2})$, which is better than in FEM. This better behavior of IGA with respect to FEM was also observed in \cite{chanMultipatchDiscontinuousGalerkin2018} for the wave equation. The plot on the right of Fig.~\ref{fig:cfl_vs_p_and_h} shows the maximum value of the time step for different mesh sizes and for degrees $p=2,3,4$. The behavior is linear with respect to the mesh size for both methods, as is expected for the employed leapfrog integrator, with some advantage to the first scheme.

To analyse the computational complexity, we compute the average wall-time per time step for the two methods considering different degrees and mesh sizes. While the choice of the linear system solver is indifferent in terms of the accuracy of the method, it is very important in terms of performance. For this reason, for the solution of the linear system associated to the pairing matrix we present results using the direct solver from MATLAB (computed with the left division operator), and exploiting the Kronecker product structure as explained in Section~\ref{sec:tensorized}. We performed the tests in the same domains as above, namely the unit cube and the coaxial cable. While in the unit cube the tensorization could be also exploited for the mass matrix, it is not possible to do it in the coaxial cable. Moreover, to further reduce the orthogonal directions, the coordinates of an internal control point of the NURBS parametrization have been slightly perturbed.
In general we have used a simulation with final time $T=2$, except for the case in which we have to perform the LU decomposition of pairing matrix directly for $p=4$ and the most refined mesh where we used a much smaller $T=0.2$ due to very long simulation times (see Remark \ref{remark:walltime}).

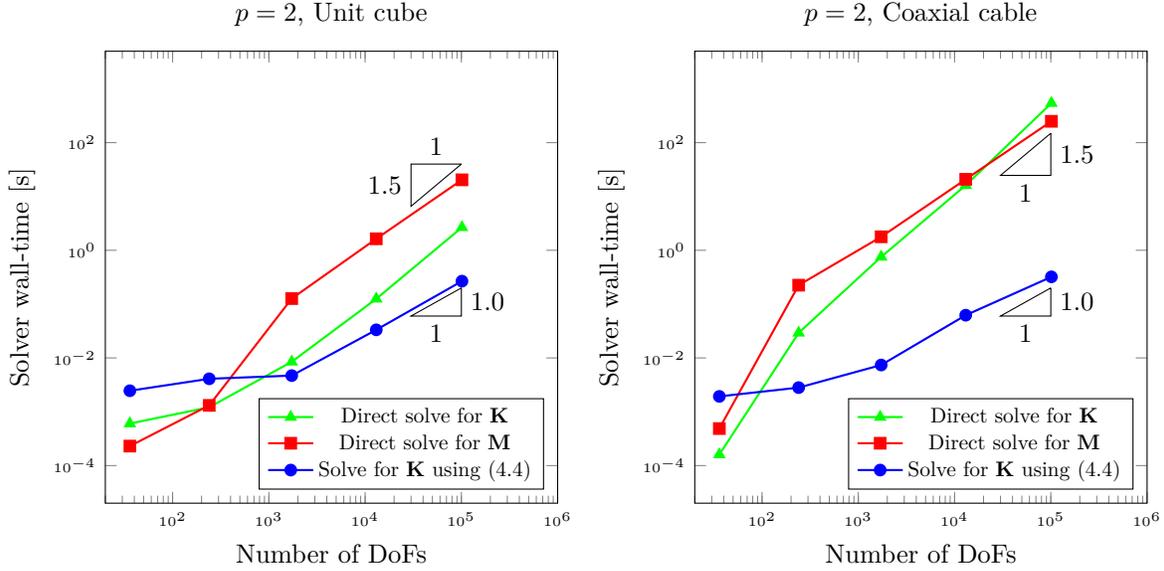
\begin{figure}[t!]
  \centering
  \begin{tikzpicture}[]
    \begin{loglogaxis}[width=0.46\textwidth, height=0.46\textwidth, 
      xlabel={Number of DoFs}, ylabel={Solver wall-time [s]}, title={$p = 2$, Unit cube},
      ymin=2E-5,ymax=5E3,
      xmin=20,xmax=1E6,
      ylabel near ticks,tick label style={font=\tiny},  legend style={font=\normalsize,nodes={scale=0.8, transform shape}}, label style={font=\normalsize},every x tick scale label/.style={at={(1,0)},anchor=north,yshift=-5pt,inner sep=0pt},legend pos= south east]
  \addplot [green,solid,thick,mark=triangle*]    table [x index=1, y index=6, col sep=comma] {figures/walltimes_final/maxwell_3d_cube_performance_direct_K_tp_jcp_p2r1_metrics.csv};
  \addplot [red,solid,thick,mark=square*] table [x index=1, y index=6, col sep=comma] {figures/walltimes_final/maxwell_3d_cube_performance_invertMsym_jcp_p2r1_metrics.csv};
  \addplot [blue,solid,thick,mark=*]   table [x index=1, y index=6, col sep=comma] {figures/walltimes_final/maxwell_3d_cube_performance_decompose_K_tp_jcp_p2r1_metrics.csv};
  \convergenceslopeinvh{3e4}{10e4}{40}{1.5}{}{};
  \convergenceslopeh{3e4}{10e4}{0.2}{1.0}{}{};
\legend{
  Direct solve for $\mathbf{K}$ \\ 
  Direct solve for $\mathbf{M}$ \\ 
  Solve for $\mathbf{K}$ using \eqref{eq:kroninverse} \\ 
}\end{loglogaxis}
\end{tikzpicture}
\begin{tikzpicture}[]
  \begin{loglogaxis}[width=0.46\textwidth, height=0.46\textwidth, 
    xlabel={Number of DoFs}, ylabel={Solver wall-time [s]}, title={$p = 2$, Coaxial cable},
    ymin=2E-5,ymax=5E3,
    xmin=20,xmax=1E6,
    ylabel near ticks,tick label style={font=\tiny},  legend style={font=\normalsize,nodes={scale=0.8, transform shape}}, label style={font=\normalsize},every x tick scale label/.style={at={(1,0)},anchor=north,yshift=-5pt,inner sep=0pt},legend pos= south east]
\addplot [green,solid,thick,mark=triangle*]    table [x index=1, y index=6, col sep=comma] {figures/walltimes_final/maxwell_3d_coaxial_performance_direct_K_tp_jcp_p2r1_metrics.csv};
\addplot [red,solid,thick,mark=square*] table [x index=1, y index=6, col sep=comma] {figures/walltimes_final/maxwell_3d_coaxial_performance_invertMsym_jcp_p2r1_metrics.csv};
\addplot [blue,solid,thick,mark=*]   table [x index=1, y index=6, col sep=comma] {figures/walltimes_final/maxwell_3d_coaxial_performance_decompose_K_tp_jcp_p2r1_metrics.csv};
\convergenceslopeh{3e4}{10e4}{150}{1.5}{}{};
\convergenceslopeh{3e4}{10e4}{0.2}{1.0}{}{};
\legend{
  Direct solve for $\mathbf{K}$ \\ 
  Direct solve for $\mathbf{M}$ \\ 
  Solve for $\mathbf{K}$ using \eqref{eq:kroninverse} \\ 
}
\end{loglogaxis}
\end{tikzpicture}
\caption{The average wall-time needed to perform a single time step in the case $p=2$. } \label{fig:performance_p2}
\end{figure}

\begin{figure}[t!]
  \centering
  \begin{tikzpicture}[]
    \begin{loglogaxis}[width=0.46\textwidth, height=0.46\textwidth, 
      xlabel={Number of DoFs}, ylabel={Solver wall-time [s]}, title={$p = 3$, Unit cube},
      ymin=2E-5,ymax=5E3,
      xmin=20,xmax=1E6,
      ylabel near ticks,tick label style={font=\tiny},  legend style={font=\normalsize,nodes={scale=0.8, transform shape}}, label style={font=\normalsize},every x tick scale label/.style={at={(1,0)},anchor=north,yshift=-5pt,inner sep=0pt},legend pos= south east]
  \addplot [green,solid,thick,mark=triangle*]    table [x index=1, y index=6, col sep=comma] {figures/walltimes_final/maxwell_3d_cube_performance_direct_K_tp_jcp_p3r2_metrics.csv};
  \addplot [red,solid,thick,mark=square*] table [x index=1, y index=6, col sep=comma] {figures/walltimes_final/maxwell_3d_cube_performance_invertMsym_jcp_p3r2_metrics.csv};
  \addplot [blue,solid,thick,mark=*]   table [x index=1, y index=6, col sep=comma] {figures/walltimes_final/maxwell_3d_cube_performance_decompose_K_tp_jcp_p3r2_metrics.csv};
  \convergenceslopeh{3e4}{10e4}{100}{1.5}{}{};
  \convergenceslopeh{3e4}{10e4}{0.5}{1.0}{}{};
\legend{
  Direct solve for $\mathbf{K}$ \\ 
  Direct solve for $\mathbf{M}$ \\ 
  Solve for $\mathbf{K}$ using \eqref{eq:kroninverse} \\ 
}
  \end{loglogaxis}
\end{tikzpicture}
\begin{tikzpicture}[]
  \begin{loglogaxis}[width=0.46\textwidth, height=0.46\textwidth, 
    xlabel={Number of DoFs}, ylabel={Solver wall-time [s]}, title={$p = 3$, Coaxial cable},
    ymin=2E-5,ymax=5E3,
    xmin=20,xmax=1E6,
    ylabel near ticks,tick label style={font=\tiny},  legend style={font=\normalsize,nodes={scale=0.8, transform shape}}, label style={font=\normalsize},every x tick scale label/.style={at={(1,0)},anchor=north,yshift=-5pt,inner sep=0pt},legend pos= south east]
\addplot [green,solid,thick,mark=triangle*]    table [x index=1, y index=6, col sep=comma] {figures/walltimes_final/maxwell_3d_coaxial_performance_direct_K_tp_jcp_p3r2_metrics.csv};
\addplot [red,solid,thick,mark=square*] table [x index=1, y index=6, col sep=comma] {figures/walltimes_final/maxwell_3d_coaxial_performance_invertMsym_jcp_p3r2_metrics.csv};
\addplot [blue,solid,thick,mark=*]   table [x index=1, y index=6, col sep=comma] {figures/walltimes_final/maxwell_3d_coaxial_performance_decompose_K_tp_jcp_p3r2_metrics.csv};
\convergenceslopeh{3e4}{10e4}{280}{1.5}{}{};
\convergenceslopeh{3e4}{10e4}{0.7}{1.0}{}{};
\legend{
  Direct solve for $\mathbf{K}$ \\ 
  Direct solve for $\mathbf{M}$ \\ 
  Solve for $\mathbf{K}$ using \eqref{eq:kroninverse} \\ 
}
\end{loglogaxis}
\end{tikzpicture}
\caption{The average wall-time needed to perform a single time step in the case $p=3$. } \label{fig:performance_p3}
\end{figure}
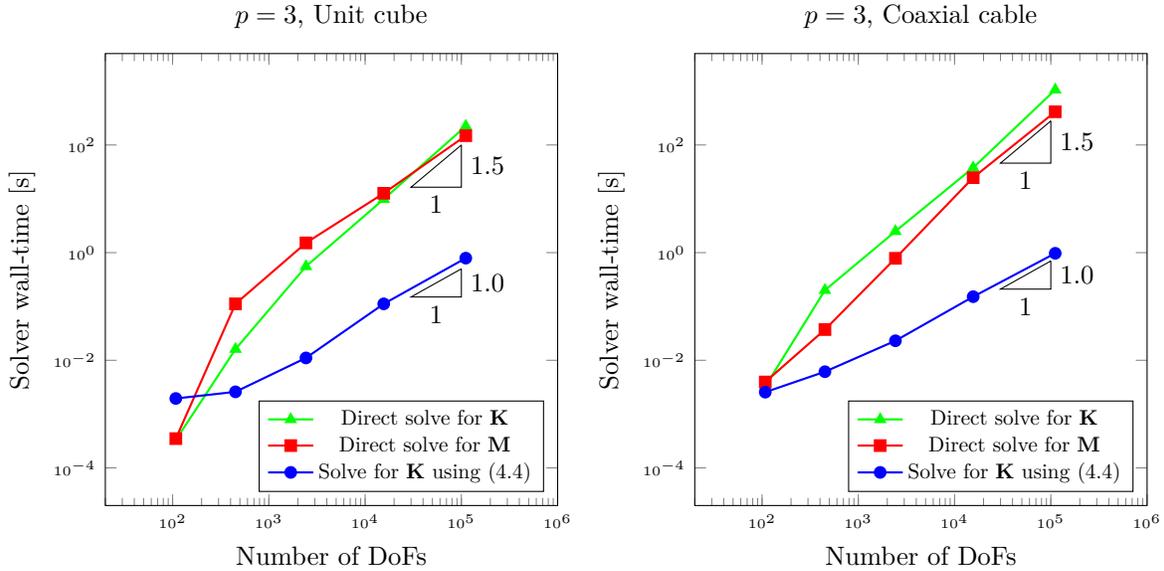

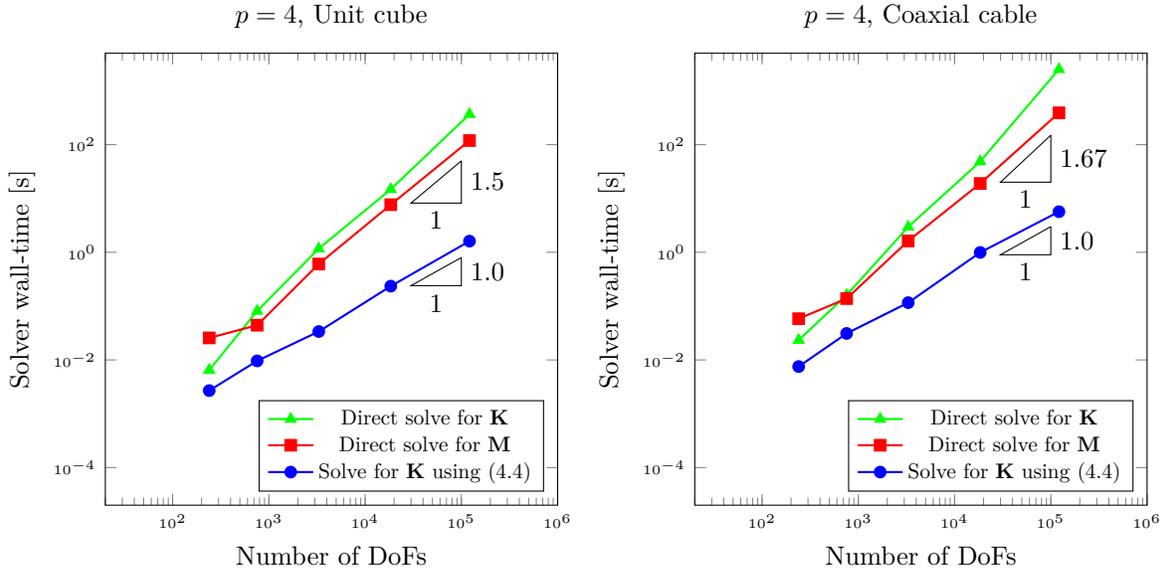
\begin{figure}[t!]
  \centering
  \begin{tikzpicture}[]
    \begin{loglogaxis}[width=0.46\textwidth, height=0.46\textwidth, 
      xlabel={Number of DoFs}, ylabel={Solver wall-time [s]}, title={$p = 4$, Unit cube},
      ymin=2E-5,ymax=5E3,
      xmin=20,xmax=1E6,
      ylabel near ticks,tick label style={font=\tiny},  legend style={font=\normalsize,nodes={scale=0.8, transform shape}}, label style={font=\normalsize},every x tick scale label/.style={at={(1,0)},anchor=north,yshift=-5pt,inner sep=0pt},legend pos= south east]
  \addplot [green,solid,thick,mark=triangle*]    table [x index=1, y index=6, col sep=comma] {figures/walltimes_final/maxwell_3d_cube_performance_direct_K_tp_jcp_p4r3_metrics.csv};
  \addplot [red,solid,thick,mark=square*] table [x index=1, y index=6, col sep=comma] {figures/walltimes_final/maxwell_3d_cube_performance_invertMsym_jcp_p4r3_metrics.csv};
  \addplot [blue,solid,thick,mark=*]   table [x index=1, y index=6, col sep=comma] {figures/walltimes_final/maxwell_3d_cube_performance_decompose_K_tp_jcp_p4r3_metrics.csv};
  \convergenceslopeh{3e4}{10e4}{50}{1.5}{}{};
  \convergenceslopeh{3e4}{10e4}{8e-1}{1.0}{}{};
\legend{
  Direct solve for $\mathbf{K}$ \\ 
  Direct solve for $\mathbf{M}$ \\ 
  Solve for $\mathbf{K}$ using \eqref{eq:kroninverse} \\ 
}
  \end{loglogaxis}
\end{tikzpicture}
\begin{tikzpicture}[]
  \begin{loglogaxis}[width=0.46\textwidth, height=0.46\textwidth, 
    xlabel={Number of DoFs}, ylabel={Solver wall-time [s]}, title={$p = 4$, Coaxial cable},
    ymin=2E-5,ymax=5E3,
    xmin=20,xmax=1E6,
    ylabel near ticks,tick label style={font=\tiny},  legend style={font=\normalsize,nodes={scale=0.8, transform shape}}, label style={font=\normalsize},every x tick scale label/.style={at={(1,0)},anchor=north,yshift=-5pt,inner sep=0pt},legend pos= south east]
    \addplot [green,solid,thick,mark=triangle*] table [x index=1, y index=6, col sep=comma] {figures/walltimes_final/maxwell_3d_coaxial_performance_invertK_jcp_p4r3_metrics.csv};
\addplot [red,solid,thick,mark=square*] table [x index=1, y index=6, col sep=comma] {figures/walltimes_final/maxwell_3d_coaxial_performance_invertMsym_jcp_p4r3_metrics.csv};
\addplot [blue,solid,thick,mark=*]   table [x index=1, y index=6, col sep=comma] {figures/walltimes_final/maxwell_3d_coaxial_performance_decompose_K_tp_jcp_p4r3_metrics.csv};
\convergenceslopeh{3e4}{10e4}{150}{1.67}{}{};
\convergenceslopeh{3e4}{10e4}{3}{1.0}{}{};
\legend{
              Direct solve for $\mathbf{K}$ \\ 
              Direct solve for $\mathbf{M}$ \\ 
              Solve for $\mathbf{K}$ using \eqref{eq:kroninverse} \\ 
}
\end{loglogaxis}
\end{tikzpicture}
\caption{The average wall-time needed to perform a single time step in the case $p=4$.} \label{fig:performance_p4}
\end{figure}

The solver wall-times are presented in Figs.~\ref{fig:performance_p2}, \ref{fig:performance_p3} and \ref{fig:performance_p4} for degrees two, three and four. 
We observe that for the new proposed approach the time stepping cost grows linearly with respect to the number of degrees of freedom (DoFs). This is consistent with the estimates of the previous section, since in three dimensions the number of DoFs scales as $\mathcal{O}(n^3)$. Instead, the computational cost that we obtain using the direct solver from MATLAB is in general better than our estimates, either when solving for the mass matrix as in the first scheme \eqref{eq:discrete1_ampere}--\eqref{eq:discrete1_hodge_mu} or for the pairing matrices as in the second scheme \eqref{eq:discrete2_ampere}--\eqref{eq:discrete2_hodge_mu}, as we would expect a quadratic growth in terms of DoFs. This is probably due to the fact that the geometry is simple, and the factorization algorithms in the left division operator inside MATLAB are able to exploit the separability of orthogonal directions. In spite of this, the new approach which exploits the tensorization is faster by roughly two orders of magnitude, and this largely compensates the larger time step due to a worse CFL condition. {Moreover, the gain would be even higher for finer meshes, because the ratio of the CFL condition for the two methods remains unchanged, while the ratio of the computational time per time step increases when refining.}

{
\begin{remark}\label{remark:walltime}
As an example, the results in the coaxial cable for $p=4$ and the finest mesh give a computational time per time step around 80 times lower for the new method compared to the standard method based on the mass matrices. Combined with the results for the CFL condition in Fig.~\ref{fig:cfl_vs_p_and_h}, a rough estimate gives a computational time around 40 times lower for the new method. This gain can be better understood comparing real numbers: the new method would run for one hour against almost two days, or e.g. one week against nine months. 
While acknowledging that these estimates are very rough, completely neglecting many aspects of practical simulations (memory access, parallelization, et cetera), they show the potential benefits of the new method.
\end{remark}
}

\section{Conclusions}\label{sec:fin}
We have developed and tested a new method for the solution of Maxwell equations with high-degree splines. The framework is based on the discretization of two de Rham complexes of differential forms and a special construction of discrete Hodge star operators between them. These Hodge operators require the solution of a linear system for a pairing matrix, which has Kronecker tensor-product structure. Exploiting this structure reduces the computational cost per time step by two orders of magnitude when compared to the solution of a linear system for the mass matrix, while the accuracy of the discretization is maintained. The method also preserves the geometry of Maxwell's equations, and in particular it conserves charges and energy. Future work directions will focus on pairing the scheme with higher order symplectic time integrators and on finding ways to exploit the inversion of the Kronecker product in \eqref{eq:system_kronecker} in a multipatch setting.

\section*{Acknowledgements}
The authors would like to thank Dr.~Mattia Tani for useful discussions regarding the solution of systems with Kronecker type matrix.
The authors acknowledge support by the Swiss National Science Foundation via the project HOGAEMS n.200021\_188589.

\bibliography{biblio}

\end{document}